\documentclass{article}
\usepackage{epsfig} % for postscript graphics files
\usepackage{psfrag}
\usepackage{amsmath} % assumes amsmath package installed
\usepackage{amssymb}  % assumes amsmath package installed

\usepackage{url}

\newtheorem{thm}{Theorem}
\newtheorem{prp}{Proposition}
\newtheorem{cor}{Corollary}
\newtheorem{ass}{Assumption}
\newtheorem{lem}{Lemma}

\newtheorem{rem}{Remark}
\newtheorem{defin}{Definition}

\newenvironment{pf}{\smallbreak\noindent{\it Proof. }}{\hfill$\Box$\smallbreak}

\usepackage{enumerate}
\usepackage{multirow}

\newcommand{\overbar}[1]{\mkern 1.5mu\overline{\mkern-1.5mu#1\mkern-1.5mu}\mkern 1.5mu}
\newcommand{\reals}{\mathbb{R}}

\newcommand{\hilbert}{\mathcal{H}}
\newcommand{\id}{{\rm{Id}}}
\DeclareMathOperator*{\argmin}{argmin}

\usepackage{tikz}
\usepackage{pgfplots}

\usetikzlibrary{calc}
\usetikzlibrary{decorations.pathreplacing}

\newcommand{\marker}[2]% [centerX,centerY,radius,angle1,angle2]
{
  \pgfgettransformentries{\myxscale}{\@tempa}{\@tempa}{\myyscale}{\@tempa}{\@tempa}
  \draw[thick] ($#1+0.08*(1/\myxscale,1/\myyscale)$)--($#1-0.08*(1/\myxscale,1/\myyscale)$);
  \draw[thick] ($#1+0.08*(-1/\myxscale,1/\myyscale)$)--($#1-0.08*(-1/\myxscale,1/\myyscale)$);
}

\newcommand{\radius}[5]% [radius,angle1,angle2]
{   
  \draw[#5] (#1,#2)--({#1+#3*cos(#4)},{#2+#3*sin(#4)});
}

\newcommand{\secant}[3]% [radius,angle1,angle2]
{   
  \draw ({#1*cos(#2)},{#1*sin(#2)})--({#1*cos(#3)},{#1*sin(#3)});
}

\newcommand{\secantcontr}[4]% [radius,angle1,angle2]
{   
  \draw ({#1*cos(#2)},{#1*sin(#2)})--({#4*#1*cos(#3)},{#4*#1*sin(#3)});
}

\newcommand{\secantalphapoint}[4]% [radius,angle1,angle,alpha]
{ 
  \coordinate (P) at
  ({(1-#4)*#1*cos(#2)+#4*#1*cos(#3)},{(1-#4)*#1*sin(#2)+#4*#1*sin(#3)});
  \marker{(P)}{0.02};
}

\newcommand{\secantalphapointtext}[6]% [radius,angle1,angle,alpha,position,text]
{ 
  \node[#5] at
  ({(1-#4)*#1*cos(#2)+#4*#1*cos(#3)},{(1-#4)*#1*sin(#2)+#4*#1*sin(#3)}) {#6};
}

\newcommand{\secantalphapointcontr}[5]% [radius,angle1,angle,alpha]
{ 
  \coordinate (P) at ({#5*((1-#4)*#1*cos(#2)+#4*#1*cos(#3))},{#5*((1-#4)*#1*sin(#2)+#4*#1*sin(#3))});
  \marker{(P)}{0.02};
}

\newcommand{\secantalphapointtextcontr}[7]% [radius,angle1,angle,alpha,position,text]
{ 
  \node[#5] at
  ({#7*((1-#4)*#1*cos(#2)+#4*#1*cos(#3))},{#7*((1-#4)*#1*sin(#2)+#4*#1*sin(#3))}) {#6};
}

\newcommand{\arccenter}[6]% [centerX,centerY,radius,angle1,angle2]
{
  \draw[#6] ([shift=(#4:#3)]#1,#2) arc (#4:#5:#3);
}

\title{\LARGE \bf
  Linear Convergence and Metric Selection for Douglas-Rachford Splitting and ADMM}

\author{Pontus Giselsson\thanks{Department of Automatic Control, Lund
    University. Email: {\tt{pontusg@control.lth.se}}.}{\phantom{a}}and Stephen
  Boyd\thanks{ Electrical Engineering Department, Stanford
    University.
    Email: {\tt{boyd@stanford.edu}}.}}

\begin{document}

\maketitle
\thispagestyle{empty}
\pagestyle{empty}

%%%%%%%%%%%%%%%%%%%%%%%%%%%%%%%%%%%%%%%%%%%%%%%%%%%%%%%%%%%%%%%%%%%%%%%%%%%%%%%% 
\begin{abstract}

  Recently, several convergence rate results for Douglas-Rachford
  splitting and the alternating direction method of multipliers (ADMM)
  have been presented in the literature. In this paper, we show global linear
  convergence rate bounds for Douglas-Rachford splitting and ADMM under strong
  convexity and smoothness assumptions. We further show that the rate bounds are tight for the class of problems
  under consideration for all feasible algorithm parameters.
  For problems that satisfy the assumptions, we show how to select step-size and metric for the algorithm that optimize the derived convergence rate bounds. For problems with a similar structure that do not satisfy the assumptions, we present heuristic step-size and metric selection methods.

\end{abstract}

%\begin{IEEEkeywords}
%Optimization algorithms, Douglas-Rachford splitting, ADMM, linear convergence.
%\end{IEEEkeywords}

%%%%%%%%%%%%%%%%%%%%%%%%%%%%%%%%%%%%%%%%%%%%%%%%%%%%%%%%%%%%%%%%%%%%%%%%%%%%%%%% 

\section{Introduction}

Optimization problems of the form
\begin{equation}
  \begin{tabular}{ll}
    minimize & $f(x)+g(x)$
  \end{tabular}
  \label{eq:splitProb}
\end{equation}
where $x$ is the variable and $f$ and $g$ are proper closed and convex,
arise in numerous applications
ranging from compressed sensing \cite{CandesCS} and statistical
estimation \cite{hastie09} to
model predictive control \cite{RawlingsBook} and medical imaging
\cite{Lustig_SparseMRI_2007}.

There exist a variety of operator splitting methods for solving problems of the form \eqref{eq:splitProb}, see \cite{Combettes2011,ParikhPA}. In this 
paper, we focus on Douglas-Rachford splitting \cite{DouglasRachford,PeacemanRachford}. Douglas-Rachford splitting is given by the iteration
\begin{align}
z^{k+1} = ((1-\alpha)\id+\alpha R_{\gamma f}R_{\gamma g})z^k,
\label{eq:DR_iter}
\end{align}
where $\alpha$ in general is constrained to satisfy $\alpha\in(0,1)$.
Here $\id$ is the identity operator, $R_{\gamma f}$ is the reflected proximal operator $R_{\gamma f}=2{\rm{prox}}_{\gamma f}-\id$, and ${\rm{prox}}_{\gamma f}$ is the proximal operator defined by
\begin{align}
{\rm{prox}}_{\gamma f}(z):=\argmin_{x}\{\gamma f(x)+\tfrac{1}{2}\|x-z\|^2\}.
\label{eq:prox_def}
\end{align}
The variable $z^k$ in \eqref{eq:DR_iter} converges to a fixed-point of $R_{\gamma f}R_{\gamma g}$ and the sequence $x^{k}={\rm{prox}}_{\gamma g}z^k$ converges to a solution of \eqref{eq:splitProb} (if it exists), see \cite[Proposition~25.6]{bauschkeCVXanal}.

When \eqref{eq:DR_iter} is applied to the Fenchel dual problem of \eqref{eq:splitProb}, this algorithm is equivalent to ADMM, see \cite{Glowinski1975,Gabay1976,BoydDistributed} for ADMM and \cite{Gabay83,EcksteinPhD} for the connection to Douglas-Rachford splitting. These methods have long been known to converge for any $\alpha\in(0,1)$ and $\gamma >0$ under mild assumptions,
\cite{Gabay1976,LionsMercier1979,EcksteinPhD}. However,
the rate of convergence in
the general case has just recently been shown to be $O(1/k)$,
\cite{DR_one_over_k_2012,even_lin_conv_2013,conv_split_schemes_Davis_2014}.

For a restricted class of monotone inclusion problems, Lions and Mercier showed in
\cite{LionsMercier1979} that the Douglas-Rachford algorithm (with $\alpha=0.5$) enjoys a
linear convergence rate. To the authors' knowledge, this was the sole
linear convergence rate results for a long period of time for these methods. Recently,
however, many works have shown linear convergence rates for
Douglas-Rachford splitting and ADMM in
different settings
\cite{HesseLuke2014,Phan2014,even_lin_conv_2013,Davis_Yin_2014,linConvADMM,GhadimiADMM,
 Panos_acc_DR_2014,lin_conv_DR_mult_block_2013,exp_conv_dist_ADMM_2013,boley_2013,
Shi_lin_conv_ADMM_2014,Bauschke_lin_rate_Friedrich,arvind_ADMM,laurent_ADMM}.

The works in \cite{HesseLuke2014,even_lin_conv_2013,boley_2013,Phan2014} concern local
linear convergence under different assumptions. The works in
\cite{exp_conv_dist_ADMM_2013,Shi_lin_conv_ADMM_2014}
consider distributed formulations and \cite{lin_conv_DR_mult_block_2013} considers multi-block ADMM,
while the works in
\cite{Davis_Yin_2014,linConvADMM,GhadimiADMM,Panos_acc_DR_2014,LionsMercier1979,Bauschke_lin_rate_Friedrich,arvind_ADMM,laurent_ADMM}
show global linear convergence. Of these, the work in
\cite{Bauschke_lin_rate_Friedrich} shows linear
convergence for Douglas-Rachford splitting when solving a subspace
intersection problem and the work in \cite{arvind_ADMM} (which appeared
online during the submission procedure of this paper) shows linear
convergence for equality constrained problems with upper and lower
bounds on the variables. The other works that show global linear convergence,
\cite{Davis_Yin_2014,linConvADMM,GhadimiADMM,Panos_acc_DR_2014,LionsMercier1979,laurent_ADMM},
show this for Douglas-Rachford splitting or ADMM under strong convexity and smoothness assumptions. 

In this paper, we provide explicit linear convergence rate bounds for the case where $f$ in \eqref{eq:splitProb} is strongly convex and smooth and we show how to select algorithm parameters to optimize the bounds. We also show that the bounds are tight for the class of problems under consideration for all feasible algorithm parameters $\gamma$ and $\alpha$. Our results generalize and/or improve corresponding results in 
\cite{Davis_Yin_2014,linConvADMM,GhadimiADMM,Panos_acc_DR_2014,LionsMercier1979,laurent_ADMM}. A detailed comparison to the 
rate bounds in these works is found in
Section~\ref{sec:comparison}. 

We also show that a relaxation factor $\alpha\geq 1$ can be used in the algorithm. We provide explicit lower and upper bounds on $\alpha$ where the upper bound depends on the smoothness and strong convexity parameters of the problem. We further show that the bounds are tight, meaning that if an $\alpha$ is chosen that is outside the allowed interval, there exists a problem from the considered class that the algorithm does not solve. 

% The possibility of using $\alpha\geq 1$ has previously been reported in the literature, \cite{laurent_ADMM}.  There, a small semi-definite program (SDP) is solved to assess whether ADMM converges for a specific strong convexity modulus, smoothness parameter, and algorithm parameters $\alpha$ and $\gamma$, as well as a guess for the linear convergence rate bound. If the SDP is feasible, the algorithm is guaranteed to convergence with that rate. They numerically note that $\alpha>1$ can be chosen, but no explicit bound is provided.

The provided convergence rate bounds for Douglas-Rachford splitting and ADMM depend on the smoothness and strong convexity parameters of the problem to be solved. These parameters depend on the space on which the problem is defined. This space can therefore be chosen by optimizing the convergence rate bound. In this paper, we show how to do this in Euclidean settings. We select a space by choosing a (metric) matrix $M$ that defines an inner-product $\langle x,y\rangle_M=x^TMy$ and its induced norm $\|x\|_M=\sqrt{\langle x,x\rangle}_M$. This matrix $M$ is chosen to optimize the convergence rate bound (typically subject to $M$ being diagonal). Letting the problem and algorithm be defined on this space can considerably improve theoretical and practical convergence. 

The metric selection can be interpreted as a preconditioning with the objective to reduce the number of iterations. A similar preconditioning is presented in \cite{GhadimiADMM} for ADMM applied to solve quadratic programs with linear inequality constraints. Our results generalize these in several directions, see Section~\ref{sec:comparison} for a detailed comparison. 

% Similar preconditioning techniques are proposed for many other iterative methods. Some focus on reducing the compuational cost within one iteration. Others (like the one we propose) focus on reducing the number of iterations in the algorithm (without increasing the per iteration cost). 
% Preconditioning for the conjugate gradient method \cite{Hestenes_CG} is primarily aimed at the latter of these objectives, see the survey \cite{Benzi_precond}. In \cite{Bramble_Uzawa,Hu_nonlin_Uzawa} both  objectives, i.e., to get cheaper and fewer iterations, are aimed at when solving linear saddle-point problems and generalizations of these. 
% For ADMM and Douglas-Rachford splitting, the former of the objectives is targeted in, e.g., \cite{ChambollePock,PockChambolle_Diag,Bredies_Sun_2015}, while the latter of the objectives is aimed at in \cite{GhadimiADMM} for ADMM and in \cite{gisBoydAut2014metric_select} for (fast) dual forward-backward splitting. 

Real-world problems rarely have the properties needed to ensure
linear convergence of Douglas-Rachford splitting or
ADMM. Therefore, we provide heuristic metric and parameter selection
methods for such problems. The heuristics cover
many problems that arise, e.g., in model
predictive control \cite{RawlingsBook}, statistical estimation
\cite{hastie09,lasso_1996}, 
compressed sensing \cite{CandesCS}, and medical
imaging \cite{Lustig_SparseMRI_2007}. 
We provide two numerical examples that show the efficiency of the theoretically justified and heuristic parameter and metric selection methods.

% A
% numerical example on a model predictive control problem is 
% provided that shows the efficiency of the proposed metric selection
% heuristic. For the
% considered problem, the execution time is decreased with about one
% order of magnitude compared to when applying the algorithm on the
% Euclidean space with the standard inner product and induced norm.

This paper extends and generalizes the conference papers
\cite{gisBoyd2014CDCprecondADMM,gisCDC2015}.

% This paper is an extended version of \cite{gisBoyd2014CDCprecondADMM}.
% It extends the results 
% in \cite{gisBoyd2014CDCprecondADMM} in that we provide convergence
% rate results for general
% Hilbert spaces, as opposed to Euclidean spaces in
% \cite{gisBoyd2014CDCprecondADMM}. Also, we
% provide a much more detailed analysis that sheds light on why our
% linear convergence rates are better than most other rates available in
% the literature.

\subsection{Notation}

We denote by $\reals$ the set of real numbers, $\reals^n$ the
set of real column-vectors of length $n$, and $\reals^{m\times n}$
the set of real matrices with $m$ rows and $n$ columns. Further
$\overbar{\reals}:=\reals\cup\{\infty\}$ denotes the extended
real line. 
Throughout this paper $\hilbert$, $\hilbert_1$, $\hilbert_2$, and $\mathcal{K}$ denote real Hilbert spaces. Their inner products are denoted by $\langle\cdot,\cdot\rangle$, the
induced norm by $\|\cdot\|$, and the identity operator by $\id$. We
specifically consider
finite-dimensional Hilbert-spaces $\mathbb{H}_H$ with inner product
$\langle x,y\rangle=x^THy$ and induced norm $\|x\|=\sqrt{x^THx}$.
We denote these by 
$\langle\cdot,\cdot\rangle_H$ and $\|\cdot\|_H$. We also denote the Euclidean
inner-product by
$\langle x,y\rangle_2=x^Ty$ and the induced norm by $\|\cdot\|_2$. The
conjugate function is denoted and defined by 
$f^{*}(y)\triangleq \sup_{x}\left\{\langle y,x\rangle-f(x)\right\}$ and the adjoint operator to a linear operator $\mathcal{A}~:~\mathcal{H}\to\mathcal{K}$ is defined as the unique operator $\mathcal{A}^*~:~\mathcal{K}\to\mathcal{H}$ that satisfies $\langle \mathcal{A}x,y\rangle=\langle x,\mathcal{A}^*y\rangle$. The range and kernel of $\mathcal{A}$ are denoted by ${\rm{ran}}\mathcal{A}$ and ${\rm{ker}}\mathcal{A}$ respectively. The orthogonal complement of a subspace $\mathcal{X}$ is denoted by $\mathcal{X}^\perp$.
The power set of a set $\mathcal{X}$, i.e., the set of all subsets of
$\mathcal{X}$, is denoted by $2^{\mathcal{X}}$. 
The graph of an
(set-valued) operator $A~:~\mathcal{X}\to 2^{\mathcal{Y}}$ is defined and denoted
by ${\rm{gph}}A = \{(x,y)\in\mathcal{X}\times\mathcal{Y}~|~y\in Ax\}$.
The inverse operator $A^{-1}$ is defined through its graph by
${\rm{gph}}A^{-1} = \{(y,x)\in\mathcal{Y}\times\mathcal{X}~|~y\in
Ax\}$. The set of fixed-points to an operator $T~:~\hilbert\to\hilbert$ is denoted and defined as ${\rm{fix}}T=\{x\in\hilbert~|~Tx=x\}$. Finally, the
class of closed, proper, and convex
functions $f~:~\hilbert\to\overbar{\reals}$ is denoted by
$\Gamma_0(\hilbert)$.

\section{Background}

In this section, we introduce some standard definitions 
that can be found, e.g. in
\cite{bauschkeCVXanal,Rockafellar1997}.

\begin{defin}[Strong monotonicity]
  An operator $A~:~\hilbert\to 2^{\hilbert}$ is
  $\sigma$-\emph{strongly monotone} with $\sigma>0$ if
  \begin{align*}
    \langle u-v,x-y\rangle\geq \sigma\|x-y\|^2
  \end{align*}
  for all $(x,u)\in{\rm{gph}}(A)$ and $(y,v)\in{\rm{gph}}(A)$.
  \label{def:str_monotone}
\end{defin}

The definition of \emph{monotonicity} is obtained by setting
$\sigma=0$ in the above definition.
% \begin{defin}[Maximal monotonicity]
%   A monotone operator $A~:~\hilbert\to 2^{\hilbert}$
%   is \emph{maximal monotone} if ${\rm{gph}}(A)$ is not a proper subset
%   of the graph of any other monotone operator
%   $B~:~\hilbert\to 2^{\hilbert}$.
% \end{defin}

In the following definitions, we
suppose that $\mathcal{D}\subseteq\hilbert$ is a nonempty subset of $\hilbert$.

\begin{defin}[Lipschitz mappings]
  A mapping $T~:~\mathcal{D}\to\hilbert$ is
  $\beta$-\emph{Lipschitz continuous} with $\beta\geq 0$ if
  \begin{align*}
    \|Tx-Ty\|\leq \beta\|x-y\|
  \end{align*}
  holds for all $x,y\in\mathcal{D}$. 
  If $\beta=1$ then $T$ is \emph{nonexpansive} and if $\beta\in[0,1)$
  then $T$ is
  $\beta$-\emph{contractive}. 
  \label{def:Lipschitz}
\end{defin}

\begin{defin}[Averaged mappings]
  A mapping $T~:~\mathcal{D}\to\hilbert$ is
  $\alpha$-\emph{averaged} if there exists a nonexpansive mapping
  $S~:~\mathcal{D}\to\hilbert$
  and $\alpha\in(0,1)$ such that $T=(1-\alpha)\id+\alpha S$.
\end{defin}

\begin{defin}[Cocoercivity]
  A mapping $T~:~\mathcal{D}\to\hilbert$ is $\beta$-cocoercive with
  $\beta > 0$ if $\beta T$ is $\tfrac{1}{2}$-averaged.
  \label{def:cocoercive}
\end{defin}

This definition implies that cocoercive mappings $T$ can be expressed as
\begin{align}
T=\tfrac{1}{2\beta}(\id+S)
\label{eq:coco_expression}
\end{align}
for some nonexpansive operator $S$.

%\subsection{Function definitions and properties}

\begin{defin}[Strong convexity]
  \label{def:strConv}
  A function $f\in\Gamma_0(\hilbert)$ is $\sigma$-\emph{strongly
    convex} with $\sigma>0$ if $f-\tfrac{\sigma}{2}\|\cdot\|^2$ is convex.
\end{defin}

A strongly convex function has a minimum curvature that is at least
$\sigma$. If $f$ is differentiable, strong convexity can equivalently
be defined as that
  \begin{align}
    f(x)\geq f(y)+\langle \nabla f(y),x-y\rangle+\tfrac{\sigma}{2}\|x-y\|^2
    \label{eq:str_cvx_def}
  \end{align}
holds for all $x,y\in\hilbert$.
Functions with a maximal curvature are called smooth. Next,
we present a smoothness definition for convex functions.

\begin{defin}[Smoothness for convex functions]
  A function $f\in\Gamma_0(\hilbert)$ is
  $\beta$-smooth with $\beta\geq 0$ if it is differentiable and
  $\tfrac{\beta}{2}\|\cdot\|^2-f$ is convex, or equivalently that 
  \begin{align}
    f(x)\leq f(y)+\langle \nabla f(y),x-y\rangle+\tfrac{\beta}{2}\|x-y\|^2
    \label{eq:fSmooth_cvx}
  \end{align}
  holds for all $x,y\in\hilbert$.
\label{def:smoothness}
\end{defin}
\begin{rem}
It can be seen from \eqref{eq:str_cvx_def} and \eqref{eq:fSmooth_cvx}
that for a function that is $\sigma$-strongly convex 
and $\beta$-smooth, we always have $\beta\geq\sigma$.
\end{rem}

\section{Douglas-Rachford splitting}

The Douglas-Rachford algorithm can be applied to solve
composite convex optimization problems of the form
\begin{align}
\begin{tabular}{ll}
minimize & $f(x)+g(x)$
\end{tabular}
\label{eq:comp_cvx}
\end{align}
where $f,g\in\Gamma_0(\hilbert)$. 
The solutions to
\eqref{eq:comp_cvx} are characterized by the following optimality
conditions, \cite[Proposition 25.1]{bauschkeCVXanal}
\begin{align}
  z&=R_{\gamma g}R_{\gamma f}z,& x&={\rm{prox}}_{\gamma f}(z)
\label{eq:DR_optcond}
\end{align}
where $\gamma>0$, the prox operator ${\rm{prox}}_{\gamma f}$ is defined in \eqref{eq:prox_def}, and
the reflected proximal operator $R_{\gamma f}=2{\rm{prox}}_{\gamma f}-\id$.
In other words, a solution to \eqref{eq:comp_cvx} is found by applying
the proximal operator on $z$, where $z$ is a fixed-point to $R_{\gamma
  g}R_{\gamma f}$.

One approach to find a fixed-point to $R_{\gamma
  g}R_{\gamma f}$ is to iterate the composition
\begin{align*}
  z^{k+1} = R_{\gamma g}R_{\gamma f}z^k.
\end{align*}
This algorithm is sometimes referred to as Peaceman-Rachford splitting,
\cite{PeacemanRachford}. However, since $R_{\gamma f}$
and $R_{\gamma g}$ are in general nonexpansive, so is their composition, and convergence of this algorithm cannot be guaranteed in the general
case. 

The Douglas-Rachford splitting algorithm is obtained by iterating
the averaged map of the nonexpansive Peaceman-Rachford operator
$R_{\gamma g}R_{\gamma f}$.
That is, it is given by the iteration
\begin{align}
  z^{k+1} = ((1-\alpha)\id+\alpha R_{\gamma g}R_{\gamma f})z^k
  \label{eq:genDougRach}
\end{align}
where $\alpha\in(0,1)$ to guarantee convergence in the general case, see \cite{Eckstein_DR_PPA}. (We will, however, see that when additional regularity assumptions are
introduced, $\alpha=1$, i.e. Peaceman-Rachford
splitting, and even some $\alpha>1$ can be used and convergence can still
be guaranteed.) 

The Douglas-Rachford algorithm \eqref{eq:genDougRach} can more
explicitly be written as
\begin{align}
\label{eq:DR1}  x^k &= {\rm{prox}}_{\gamma f}(z^k)\\
\label{eq:DR2}  y^k &= {\rm{prox}}_{\gamma g}(2x^k-z^k)\\
\label{eq:DR3}  z^{k+1} &= z^k+2\alpha(y^k-x^k)
\end{align}
Note that sometimes, the algorithm obtained by letting $\alpha=\tfrac{1}{2}$ is called Douglas-Rachford splitting \cite{EcksteinPhD}. Here we use the name Douglas-Rachford splitting for all feasible values of $\alpha$.

% known as Douglas-Rachford splitting is obtained by letting
% $\alpha=\tfrac{1}{2}$ in \eqref{eq:genDougRach}, but we will use the
% term Douglas-Rachford splitting for all values of $\alpha$.

\subsection{Linear convergence}

Under some regularity assumptions, the convergence of the
Douglas-Rachford algorithm is linear. We will analyze its convergence
under the following set of assumptions:
\begin{ass}
Suppose that:
\begin{enumerate}[(i)]
\item $f$ and $g$ are proper, closed, and convex.
\item $f$ is $\sigma$-strongly convex and $\beta$-smooth.
\end{enumerate}
\label{ass:prob_ass}
\end{ass}
To show linear convergence rates of the Douglas-Rachford algorithm
under these regularity assumptions on $f$, we need to characterize some
properties of the proximal and reflected proximal operators to $f$.
Specifically, we will show that the reflected proximal operator of $f$
is contractive (as opposed to nonexpansive in the general case). We will
also provide a tight contraction factor.

The key to arriving at this contraction factor is the following, to the authors' knowledge, novel
(but straightforward) interpretation of the proximal operator.
\begin{prp}
  Assume that $f\in\Gamma_0(\hilbert)$ and that $\gamma\in(0,\infty)$ and define $f_{\gamma}$ as 
\begin{align}
f_{\gamma} := \gamma f+\tfrac{1}{2}\|\cdot\|^2.
\label{eq:reg_fcn}
\end{align}
Then ${\rm{prox}}_{\gamma f}(y)
  = \nabla f_\gamma^*(y)$.
  \label{prp:prox_to_conj}
\end{prp}
\begin{pf}
Since the proximal operator is the resolvent of $\gamma\partial f$,
see \cite[Example~23.3]{bauschkeCVXanal},
  we have ${\rm{prox}}_{\gamma f}(y) = (\id+\gamma \partial
  f)^{-1}y=(\partial f_\gamma)^{-1}y$. Since $f\in\Gamma_0(\hilbert)$ and $\gamma\in(0,\infty)$
  also $f_{\gamma}\in\Gamma_0(\hilbert)$. Therefore
  \cite[Corollary~16.24]{bauschkeCVXanal} implies that
  ${\rm{prox}}_{\gamma f}(y)=(\partial f_\gamma)^{-1}y=\nabla
  f_{\gamma}^*(y)$, where differentiability of $f_{\gamma}^*$ follows
  from \cite[Theorem~18.15]{bauschkeCVXanal}, since $f_{\gamma}$ is
  1-strongly convex, and since $f=(f^*)^*$ for proper, closed, and convex
  functions, see \cite[Theorem~13.32]{bauschkeCVXanal}. This concludes
  the proof.
\end{pf}

This interpretation of the proximal operator is used to prove the
following proposition. The proof is found in Appendix~\ref{app:prox_coco_pf}.
\begin{prp}
  Assume that $f\in\Gamma_0(\hilbert)$ is $\sigma$-strongly convex and
  $\beta$-smooth and that $\gamma\in (0,\infty)$. Then ${\rm{prox}}_{\gamma f} -\tfrac{1}{1+\gamma\beta} \id$ is
  $\tfrac{1}{\tfrac{1}{1+\gamma\sigma}-\tfrac{1}{1+\gamma\beta}}$-cocoercive if $\beta>\sigma$ and 0-Lipschitz if $\beta=\sigma$. 
\label{prp:prox_coco}
\end{prp}

This result is used to
show the following contraction properties of the reflected proximal
operator. A proof to this result, which is one of the main results of
the paper, is found in Appendix~\ref{app:refl_res_contr_pf}. 
\begin{thm}
  Suppose that $f\in\Gamma_0(\hilbert)$ is $\sigma$-strongly convex and $\beta$-smooth and that $\gamma\in(0,\infty)$. Then
  $R_{\gamma f}$ is
  $\max\left(\tfrac{\gamma \beta-1}{\gamma\beta+1},\tfrac{1-\gamma \sigma}{\gamma \sigma+1}\right)$-contractive.
  \label{thm:refl_res_contr_subdiff}
\end{thm}

For future reference, we let this contraction factor be denoted by $\delta$, i.e.,
\begin{align}
\delta := \max\left(\tfrac{\gamma \beta-1}{\gamma\beta+1},\tfrac{1-\gamma \sigma}{\gamma \sigma+1}\right).
\label{eq:delta}
\end{align}

Theorem~\ref{thm:refl_res_contr_subdiff} lays the foundation for the linear convergence rate
result in the following theorem, which is proven in
Appendix~\ref{app:DR_lin_conv_subdiff_pf}.
\begin{thm}
  Suppose that Assumption~\ref{ass:prob_ass} holds, that $\gamma\in(0,\infty)$, that $\alpha\in (0,\tfrac{2}{1+\delta})$ with $\delta$ in \eqref{eq:delta}. Then the
  Douglas Rachford algorithm \eqref{eq:genDougRach} converges linearly
  towards a fixed-point $\bar{z}\in{\rm{fix}}(R_{\gamma g}R_{\gamma f})$ at least with rate
  $|1-\alpha|+\alpha\delta$, i.e.,:
  \begin{align*}
    \|z^{k+1}-\bar{z}\|\leq\left(|1-\alpha|+\alpha\delta\right)\|z^{k}-\bar{z}\|.
  \end{align*}
  \label{thm:DR_lin_conv_subdiff}
\end{thm}
\begin{rem}
Note that $\alpha> 1$ can be
used in the Douglas-Rachford algorithm
when solving problems that satisfy Assumption~\ref{ass:prob_ass}. (This also holds for general relaxed iteration of a contractive mapping.)
That $\alpha$-values greater than 1 can be used is reported in \cite{laurent_ADMM} for ADMM (i.e., for dual Douglas-Rachford). They solve a small SDP to assess whether ADMM converges with a specific rate, for specific problem parameters $\beta$ and $\sigma$, and specific algorithm parameters $\alpha$ and $\gamma$. This SDP gives an affirmative answer for some problems and some $\alpha>1$. As opposed to here, no explicit bounds on $\alpha$ are provided.
\label{rem:alpha_upper_bound}
\end{rem}

Also the $x^k$ iterates in \eqref{eq:DR1} converge linearly. The following corollary is proven in Appendix~\ref{app:x_lin_conv_pf}.
\begin{cor}
Let $x^\star$ be the (unique) solution to \eqref{eq:comp_cvx}. Then the $x^k$ iterates in \eqref{eq:DR1} satisfy
\begin{align}
\|x^{k+1}-x^\star\|\leq(|1-\alpha|+\delta\alpha)^{k+1}\tfrac{1}{1+\gamma\sigma}\|z^0-\bar{z}\|.
\label{eq:xk_R_lin}
\end{align}
\label{cor:x_lin_conv}
\end{cor}

\begin{rem}
Note that Corollary~\ref{cor:x_lin_conv} and Theorem~\ref{thm:DR_lin_conv_subdiff} still hold if the order of $R_{\gamma g}$ and $R_{\gamma f}$ is swapped in the Douglas-Rachford algorithm \eqref{eq:genDougRach}.
\end{rem}

We can choose the algorithm parameters $\gamma$ and $\alpha$ to
optimize the bound on  
the convergence rate in Theorem~\ref{thm:DR_lin_conv_subdiff}. This is
done in the following proposition.
\begin{prp}
  Suppose that Assumption~\ref{ass:prob_ass} holds. Then the optimal
  parameters for the 
  Douglas-Rachford algorithm in \eqref{eq:genDougRach} are 
  $\alpha=1$ and $\gamma = 
  \tfrac{1}{\sqrt{\sigma\beta}}$. Further, the optimal rate is 
  $\tfrac{\sqrt{\beta/\sigma}-1}{\sqrt{\beta/\sigma}+1}$.
\label{prp:opt_param_subdiff}
\end{prp}
\begin{pf}
Since $\delta\in[0,1)$, the rate in Theorem~\ref{thm:DR_lin_conv_subdiff},
$|1-\alpha|+\alpha\delta$,
is a decreasing function of $\alpha$ for $\alpha\leq 1$ and increasing
for $\alpha \geq 1$. Therefore the rate
factor is optimized by $\alpha=1$. The
$\gamma$ parameter should be chosen to minimize the max-expression
$\max\left(\tfrac{\gamma\beta-1}{\gamma\beta+1},\tfrac{1-\gamma\sigma}{1+\gamma\sigma}\right)$ defining $\delta$ in \eqref{eq:delta}.
This is done by setting the arguments equal, which gives
$\gamma=1/\sqrt{\beta\sigma}$. Inserting these values into the rate factor
expression gives $\tfrac{\sqrt{\beta/\sigma}-1}{\sqrt{\beta/\sigma}+1}$.
\end{pf}

\begin{rem}
Note that $\alpha=1$ is optimal in Proposition~\ref{prp:opt_param_subdiff}. So
the Peaceman-Rachford algorithm gives the best bound on the
convergence rate under Assumption~\ref{ass:prob_ass},
even though it is not guaranteed to
converge in the general case. The reason is that 
$R_{\gamma f}$ becomes contractive under Assumption~\ref{ass:prob_ass} (as opposed to nonexpansive in the general case).
\end{rem}

\subsection{Tightness of bounds}
\label{sec:tightness_DR}

In this section, we show that the linear convergence rate bounds in Theorem~\ref{thm:DR_lin_conv_subdiff} are tight. We consider a two dimensional example of the form \eqref{eq:comp_cvx} in a standard Euclidean space to show tightness. Let $x=(x_1,x_2)$, then the functions used are
\begin{align}
\label{eq:f_def} f(x) &=\tfrac{1}{2}(\beta x_1^2+\sigma x_2^2),\\
\label{eq:g_def1} g_1(x) &=0,\\
\label{eq:g_def2} g_2(x) &=\iota_{x=0}(x),
\end{align}
where $\beta\geq\sigma>0$ and $\iota_{x=0}$ is the indicator function for $x=0$, i.e., $\iota_{x=0}(x)=0$ iff $x=0$, otherwise $\infty$.

The function $f$ is $\beta$-smooth since $\tfrac{\beta}{2}\|x\|_2^2-f(x)=\tfrac{\beta-\sigma}{2}x_2^2$ is convex. It is $\sigma$-strongly convex since
$f(x)-\tfrac{\sigma}{2}\|x\|_2^2=\tfrac{\beta-\sigma}{2}x_1^2$
is convex. So the function $f$ satisfies Assumption~\ref{ass:prob_ass}(ii).

The proximal operator to $f$ is
\begin{align}
\nonumber{\rm{prox}}_{\gamma f}(y) &=\argmin_{x}\{\gamma f(x)+\tfrac{1}{2}\|x-y\|_2^2\}\\
\nonumber&=\argmin_{x}\{\tfrac{1}{2}\left((\beta\gamma+1) x_1^2+(\sigma\gamma+1) x_2^2\right) + x_1y_1+x_2y_2\}\\
\label{eq:f_prox}&=(\tfrac{1}{1+\gamma\beta}y_1,\tfrac{1}{1+\gamma\sigma}y_2)
\end{align}
where $y=(y_1,y_2)$. The reflected proximal operator is
\begin{align}
\nonumber R_{\gamma f}(y) &= 2{\rm{prox}}_{\gamma f}(y)-y\\
\nonumber &=2(\tfrac{1}{1+\gamma\beta}y_1,\tfrac{1}{1+\gamma\sigma}y_2)-(y_1,y_2)\\
\label{eq:f_refl_prox}&=(\tfrac{1-\gamma\beta}{1+\gamma\beta}y_1,\tfrac{1-\gamma\sigma}{1+\gamma\sigma}y_2).
\end{align}

The proximal and reflected proximal operators to $g_1$ in \eqref{eq:g_def1} are ${\rm{prox}}_{\gamma g_1} = R_{\gamma g_1} = \id$. The proximal and reflected proximal operators to $g_2$ in \eqref{eq:g_def2} are ${\rm{prox}}_{\gamma g_2}(x) = 0$ and $R_{\gamma g_2}=2{\rm{prox}}_{\gamma g_2} -\id=-\id$.

Next, these results are used to show tightness of the
convergence rate bounds in Theorem~\ref{thm:DR_lin_conv_subdiff}. The tightness result is proven in Appendix~\ref{app:tight_bound_pf}. 

\begin{thm}
The convergence rate bound in Theorem~\ref{thm:DR_lin_conv_subdiff} for the Douglas-Rachford splitting algorithm \eqref{eq:genDougRach} is tight for the class of problems that satisfy Assumption~\ref{ass:prob_ass} for all algorithm parameters specified in Theorem~\ref{thm:DR_lin_conv_subdiff}, namely $\alpha\in(0,\tfrac{2}{1+\delta})$ with $\delta$ in \eqref{eq:delta} and $\gamma\in(0,\infty)$. If instead $\alpha\not\in(0,\tfrac{2}{1+\delta})$, there exist a problem that satisfies Assumption~\ref{ass:prob_ass} that the Douglas-Rachford algorithm does not solve.
\label{thm:tight_bound}
\end{thm}
\begin{rem}
This result on tightness can be generalized to any real Hilbert space. This is obtained by considering the same $g_1$ and $g_2$ as in \eqref{eq:g_def1} and \eqref{eq:g_def2} but with 
\begin{align*}
f(x) &=\sum_{i=1}^{|\hilbert|}\tfrac{\lambda_i}{2}\langle
x,\phi_i\rangle^2.
\end{align*}
Here $\{\phi_i\}_{i=1}^{|\hilbert|}$ is an orthonormal basis for $\hilbert$,
$|\hilbert|$ is the dimension of the space $\mathcal{H}$
(possibly infinite), and $\lambda_i=\sigma>0$
or $\lambda_i=\beta\geq\sigma$ for all $i$, where at least one $\lambda_i=\sigma$ and one $\lambda_i=\beta$.
\label{rem:tightness_hilbert}
\end{rem}

\section{ADMM}
\label{sec:ADMM}

In this section, we apply the results from the previous section to the Fenchel dual problem. Since ADMM applied to the primal problem is equivalent to the Douglas-Rachford algorithm applied to the dual problem, the results obtained in this section hold for ADMM. 

We consider solving problems of the form 
\begin{equation}
  \begin{tabular}{ll}
    minimize & $f(x)+g(y)$\\
subject to &$\mathcal{A}x+\mathcal{B}y=c$
  \end{tabular}
  \label{eq:prob_wA}
\end{equation}
where $f\in\Gamma_0(\hilbert_1)$, $g\in\Gamma_0(\hilbert_2)$, $\mathcal{A}~:~\hilbert_1\to\mathcal{K}$ and $\mathcal{B}~:~\hilbert_2\to\mathcal{K}$ are bounded linear operators and $c\in\mathcal{K}$. In addition, we assume:
\begin{ass}~
  \begin{enumerate}[(i)]
  \item $f\in\Gamma_0(\hilbert_1)$ is $\beta$-smooth and
    $\sigma$-strongly convex.
  \item $\mathcal{A}~:~\hilbert_1\to\mathcal{K}$ is surjective.
  \end{enumerate}
  \label{ass:probAss_wA}
\end{ass}
The assumption that $\mathcal{A}$ is
a surjective bounded linear operator reduces to that $\mathcal{A}$ is a
real matrix with full row rank in the Euclidean case.

Problems of the form \eqref{eq:prob_wA} cannot be directly efficiently
solved using Douglas-Rachford splitting. Therefore, we
instead solve the (negative) Fenchel dual problem, which is
\begin{align}
  \begin{tabular}{ll}
    minimize & $d_1(\mu)+d_2(\mu)$
  \end{tabular}
  \label{eq:dualProb}
\end{align}
where $d_1,d_2\in\Gamma_0(\mathcal{K})$ are
\begin{align}
  d_1(\mu)&:= f^*(-\mathcal{A}^*\mu)+\langle c,\mu\rangle,&d_2(\mu)=g^*(-\mathcal{B}^*\mu).
  \label{eq:dFcn}
\end{align}
The Douglas-Rachford algorithm when solving the dual problem becomes
\begin{align}
z^{k+1}=((1-\alpha)\id+\alpha R_{\gamma d_1}R_{\gamma d_2})z^k.
\label{eq:DR_dual}
\end{align}
This formulation will be used when analyzing ADMM since for $\alpha=\tfrac{1}{2}$ it is equivalent to ADMM and for $\alpha\in(0,\tfrac{1}{2}]$ and $\alpha\in[\tfrac{1}{2},1)$ it is equivalent to under- and over-relaxed ADMM,  see \cite{Gabay83,EcksteinPhD,Eckstein_DR_PPA}. Note that we only use \eqref{eq:DR_dual} as an analysis algorithm for ADMM. The algorithm should still be implemented as the standard ADMM algorithm (with over- or under-relaxation). Here we state ADMM with scaled dual variables $u$, see \cite{BoydDistributed}:
\begin{align}
\label{eq:ADMM1}x^{k+1} &= \argmin_{x}\{f(x)+\tfrac{\gamma}{2}\|Ax+By^{k}-c+u^k\|_2^2\}\\
\label{eq:ADMM2}x_A^{k+1}&=2\alpha Ax^{k+1}-(1-2\alpha)(By^{k}-c)\\
\label{eq:ADMM3}y^{k+1} &= \argmin_{y}\{g(y)+\tfrac{\gamma}{2}\|x_A^{k+1}+By-c+u^k\|_2^2\}\\
\label{eq:ADMM4}u^{k+1} &= u^k+ (x_A^{k+1}+By^{k+1}-c)
\end{align}
where $z^{k}$ in \eqref{eq:DR_dual} satisfies $z^{k}=\gamma(u^{k}-By^k)$, see \cite[Theorem~8]{Eckstein_DR_PPA} and \cite[Appendix~B]{gis_line_search}.

\subsection{Linear convergence}

To show linear convergence rate results in this dual setting, we need to quantify the
strong convexity and smoothness parameters for $d_1$ in \eqref{eq:dFcn}. This is done in
the following proposition.
\begin{prp}
  Suppose that Assumption~\ref{ass:probAss_wA} holds. Then
  $d_1\in\Gamma_0(\mathcal{K})$ in \eqref{eq:dFcn} is $\tfrac{\|\mathcal{A}^*\|^2}{\sigma}$-smooth and
  $\tfrac{\theta^2}{\beta}$-strongly convex, where $\theta>0$ always exists and satisfies
  $\|\mathcal{A}^*\mu\|\geq \theta\|\mu\|$ for all
  $\mu\in\mathcal{K}$. 
\label{prp:dual_fcn_prop_gen}
\end{prp}
\begin{pf}
We define $d(\mu):=f^*(-\mathcal{A}^*\mu)$, so $d_1(\mu)=d(\mu)+\langle c,\mu\rangle$. We first note that the linear term $\langle c,\mu\rangle$ does not affect the strong convexity or smoothness parameters. So, we proceed by showing $d$ is $\tfrac{\|\mathcal{A}^*\|^2}{\sigma}$-smooth and
  $\tfrac{\theta^2}{\beta}$-strongly convex.

  Since $f$ is $\sigma$-strongly convex,
  \cite[Theorem~18.15]{bauschkeCVXanal} gives that $f^*$ is
  $\tfrac{1}{\sigma}$-smooth and that $\nabla f^*$ is
  $\tfrac{1}{\sigma}$-Lipschitz continuous. Therefore, $\nabla d$ satisfies
  \begin{align*}
    \|\nabla d(\mu)-\nabla d(\nu)\| &=\|\mathcal{A}\nabla f^*(-\mathcal{A}^*\mu)-\mathcal{A}\nabla
    f^*(-\mathcal{A}^*\nu)\| \\
    &\leq \tfrac{\|\mathcal{A}\|}{\sigma}\|\mathcal{A}^*(\mu-\nu)\| 
    \\&\leq \tfrac{\|\mathcal{A}^*\|^2}{\sigma}\|\mu-\nu\|
  \end{align*}
  since $\|\mathcal{A}\|=\|\mathcal{A}^*\|$. This is equivalent to
  that $d$ is $\tfrac{\|\mathcal{A}^*\|^2}{\sigma}$-smooth, see
  \cite[Theorem~18.15]{bauschkeCVXanal}.

  Next, we show the strong convexity result for $d$. Since $f$ is
  $\beta$-smooth, $f^*$ is
  $\tfrac{1}{\beta}$-strongly convex, and $\nabla f^*$ is
  $\tfrac{1}{\beta}$-strongly monotone,  see \cite[Theorem~18.15]{bauschkeCVXanal}.
  Therefore, $\nabla d$ satisfies
  \begin{align*}
    \langle \nabla d(\mu)-\nabla d(\nu),\mu-\nu\rangle 
    &=\langle -\mathcal{A}(\nabla f^*(-\mathcal{A}^*\mu)-\nabla
    f^*(-\mathcal{A}^*\nu)),\mu-\nu\rangle \\
    &=\langle \nabla f^*(-\mathcal{A}^*\mu)-\nabla
    f^*(-\mathcal{A}^*\nu),-\mathcal{A}^*\mu+\mathcal{A}^*\nu)\rangle \\
    &\geq\tfrac{1}{\beta}
    \|\mathcal{A}^*(\mu-\nu)\|^2\\&\geq \tfrac{\theta^2}{\beta}\|\mu-\nu\|^2.
  \end{align*}
This, by \cite[Theorem~18.15]{bauschkeCVXanal}, is equivalent to $d$ being
$\tfrac{\theta^2}{\beta}$-strongly convex.
  That $\theta>0$ follows from \cite[Fact~2.18 and Fact~2.19]{bauschkeCVXanal}.
  Specifically, \cite[Fact~2.18]{bauschkeCVXanal} says that
  ${\rm{ker}}\mathcal{A}^*=({\rm{ran}}\mathcal{A})^\perp=\emptyset$, since $\mathcal{A}$ is
  surjective. Since ${\rm{ran}}\mathcal{A}=\mathcal{K}$ (again by
  surjectivity), it is closed. Then \cite[Fact~2.19]{bauschkeCVXanal} states that there
  exists $\theta>0$ such that $\|\mathcal{A}^*\mu\|\geq
  \theta\|\mu\|$ for all
  $\mu\in({\rm{ker}} \mathcal{A}^*)^\perp=(\emptyset)^\perp=\mathcal{K}$.

\end{pf}

Combining this result with Theorem~\ref{thm:refl_res_contr_subdiff} implies that $R_{\gamma d_1}$ is $\hat{\delta}$-contractive
with
\begin{align}
\hat{\delta}&:=\max\left(\tfrac{\gamma\hat{\beta}-1}{1+\gamma\hat{\beta}},\tfrac{1-\gamma\hat{\sigma}}{1+\gamma\hat{\sigma}}\right)
\label{eq:deltahat}
\end{align}
where $\hat{\beta}=\tfrac{\|\mathcal{A}^*\|^2}{\sigma}$ and $\hat{\sigma} = \tfrac{\theta^2}{\beta}$.
This gives the following immediate corollary to Theorem~\ref{thm:DR_lin_conv_subdiff} and Proposition~\ref{prp:opt_param_subdiff}.
\begin{cor}
  Suppose that Assumption~\ref{ass:probAss_wA}
  holds, that $\gamma\in(0,\infty)$, and that $\alpha\in(0,\tfrac{2}{1+\hat{\delta}})$ with $\hat{\delta}$ in \eqref{eq:deltahat}. Then algorithm \eqref{eq:DR_dual} (or equivalently ADMM applied to solve \eqref{eq:prob_wA}) converges at least
  with rate $|1-\alpha|+\alpha\hat{\delta}$, i.e.,
\begin{align*}
\|z^{k+1}-\bar{z}\|\leq\left(|1-\alpha|+\alpha\hat{\delta}\right)\|z^k-\bar{z}\|
\end{align*}
where $\bar{z}\in{\rm{fix}}(R_{\gamma d_1}R_{\gamma d_2})$.

Further, the algorithm parameters $\gamma$
  and $\alpha$ that optimize
  the rate bound are $\alpha=1$ and $\gamma = \tfrac{1}{\sqrt{\hat{\beta}\hat{\sigma}}}=
  \tfrac{\sqrt{\beta\sigma}}{\|\mathcal{A}^*\|\theta}$. The
  optimized linear convergence rate bound factor is
  $\tfrac{\sqrt{\hat{\kappa}}-1}{\sqrt{\hat{\kappa}}+1}$, where $\hat{\kappa} =
  \tfrac{\hat{\beta}}{\hat{\sigma}} = \tfrac{\|\mathcal{A}^*\|^2\beta}{\theta^2\sigma}$.
\label{cor:ADMM_lin_conv_subdiff}
\end{cor}
% \begin{pf}
%   This follows directly from Propositions~\ref{prp:dual_fcn_prop_gen}
%   and~\ref{prp:opt_param_subdiff} and
%   Theorem~\ref{thm:DR_lin_conv_subdiff}.
% \end{pf}

\begin{rem}
Similarly to in the primal setting in Corollary~\ref{cor:x_lin_conv}, R-linear convergence for the $x^k$ update in \eqref{eq:ADMM1} can be shown. The proof is more elaborate than the corresponding proof for primal Douglas-Rachford. This result is therefore omitted due to space considerations.
\end{rem}

\subsection{Tightness of bounds}
\label{sec:tigthness_ADMM}

Also in this dual setting, the convergence rate estimates are tight. Consider again the functions \eqref{eq:f_def}, \eqref{eq:g_def1}, and \eqref{eq:g_def2}. Further let $c=0$, $\mathcal{B}=-I$, and $\mathcal{A}:=\begin{bmatrix}\begin{smallmatrix} \theta & 0\\ 0 & \zeta\end{smallmatrix}\end{bmatrix}$ with $\zeta\geq\theta>0$.

Since $f^*(\nu) = \tfrac{1}{2}(\tfrac{1}{\beta}\nu_1^2+\tfrac{1}{\sigma}\nu_2^2)$ for $f$ in \eqref{eq:f_def}, we get
\begin{align*}
d_1(\mu)=f^*(-\mathcal{A}^T\mu)=\tfrac{1}{2}(\tfrac{\theta^2}{\beta}\mu_1^2+\tfrac{\zeta^2}{\sigma}\mu_2^2).
\end{align*}
Since $d_1$ is on the same form as $f$ in Section~\ref{sec:tightness_DR}, $d_1$ is $\hat{\sigma}:=\tfrac{\theta^2}{\beta}$-strongly convex and $\hat{\beta}:=\tfrac{\zeta^2}{\sigma}$-smooth, or equivalently $\hat{\beta}:=\tfrac{\|\mathcal{A}^T\|^2}{\sigma}$-smooth since $\zeta=\|\mathcal{A}\|=\|\mathcal{A}^T\|$. Further $d_2(\mu)=g^*(-\mathcal{B}^T\mu)=g^*(\mu)$ where either $g=g_1$ with $g_1$ in \eqref{eq:g_def1} or $g=g_2$ with $g_2$ in \eqref{eq:g_def2}. The functions $g_1$ and $g_2$ in \eqref{eq:g_def1} and \eqref{eq:g_def2} are each other's conjugates, i.e., $g_1^*=g_2$ and $g_2^*=g_1$. Therefore, the dual problem
\begin{align*}
\begin{tabular}{ll}
minimize & $d_1(\mu)+d_2(\mu)$
\end{tabular}
\end{align*}
(with $d_2=g_1$ or $d_2=g_2$) is the same as the problem considered in Section~\ref{sec:tightness_DR} but with the smoothness parameter $\beta$ in Section~\ref{sec:tightness_DR} replaced by $\hat{\beta}=\tfrac{\|\mathcal{A}^*\|^2}{\sigma}$ and the strong convexity modulus $\sigma$ in Section~\ref{sec:tightness_DR} replaced by $\hat{\sigma}=\tfrac{\theta^2}{\beta}$. Therefore, we can state the following immediate corollary to Theorem~\ref{thm:tight_bound}.
\begin{cor}
The convergence rate bound in Corollary~\ref{cor:ADMM_lin_conv_subdiff} for ADMM applied to the primal or, equivalently, the Douglas-Rachford algorithm applied to the dual \eqref{eq:DR_dual}, is tight for the class of problems that satisfy Assumption~\ref{ass:probAss_wA} for all algorithm parameters specified in Corollary~\ref{cor:ADMM_lin_conv_subdiff}, namely for $\alpha\in(0,\tfrac{2}{1+\hat{\delta}})$ with $\hat{\delta}$ in \eqref{eq:deltahat} and $\gamma\in(0,\infty)$. If instead $\alpha\not\in(0,\tfrac{2}{1+\hat{\delta}})$, there exist a problem that satisfies Assumption~\ref{ass:probAss_wA} that the algorithm does not solve.
\label{cor:tight_bound}
\end{cor}
\begin{rem}
As in the primal Douglas-Rachford case, tightness can be shown in general real Hilbert spaces. This is obtained with the same $f$, $g_1$, and $g_2$ as in Remark~\ref{rem:tightness_hilbert} and with $\mathcal{B}=-\id$, $c=0$, and 
\begin{align*}
 \mathcal{A}x &=\sum_{i=1}^{|\hilbert|}\nu_i\langle x,\phi_i\rangle\phi_i
\end{align*}
where $\nu_i=\theta>0$ if $\lambda_i=\beta$ in Remark~\ref{rem:tightness_hilbert} and $\nu_i=\zeta\geq\theta$ if $\lambda_i=\sigma$ in Remark~\ref{rem:tightness_hilbert}. 
\end{rem}

\section{Related work}
\label{sec:comparison}

\begin{figure}
\centering
\includegraphics{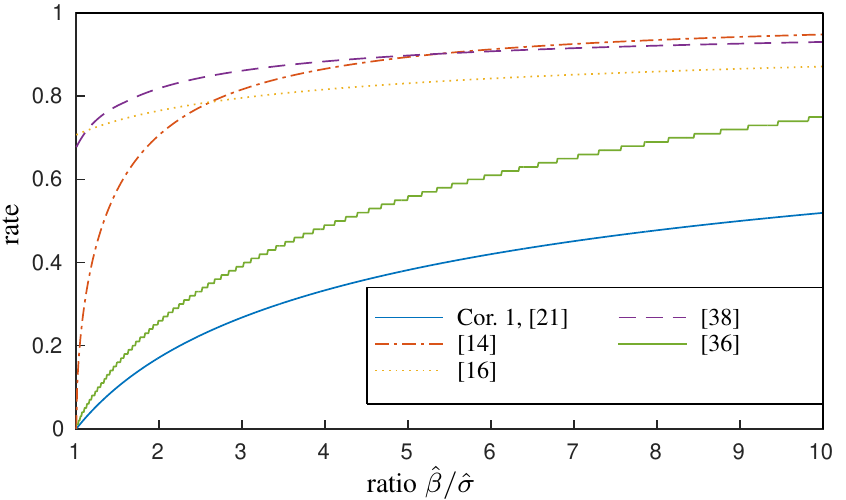}
\caption{Comparison between linear convergence rate bounds for
  Douglas-Rachford splitting/ADMM provided in
  \cite{Davis_Yin_2014,linConvADMM,GhadimiADMM,Panos_acc_DR_2014,laurent_ADMM} and in Corollary~\ref{cor:ADMM_lin_conv_subdiff}.}
\label{fig:rate_fig}
\end{figure}

Recently, many works have appeared that show linear convergence rate bounds for Douglas-Rachford splitting and ADMM under different sets of assumptions. In this section, we will discuss and compare some of these linear convergence rate bounds that hold with the assumptions stated in this paper. For simplicity, we will compare the various optimal rates, i.e., rates obtained by using optimal algorithm parameters.

Specifically, we will compare our results in 
Proposition~\ref{prp:opt_param_subdiff} and Corollary~\ref{cor:ADMM_lin_conv_subdiff} to the previously known linear convergence rate results in
\cite{Davis_Yin_2014,linConvADMM,GhadimiADMM,Panos_acc_DR_2014,LionsMercier1979}
and the linear convergence rate \cite{laurent_ADMM} that appeared
online during the submission procedure of this paper. 

Some of these results consider the Douglas-Rachford algorithm while others consider ADMM. A Douglas-Rachford rate result becomes an ADMM rate result by replacing the strong convexity and smoothness parameters $\sigma$ and $\beta$ with the dual problem counterparts $\hat{\sigma}=\tfrac{\theta^2}{\beta}$
and $\hat{\beta}=\tfrac{\|\mathcal{A}^*\|^2}{\sigma}$, see Proposition~\ref{prp:dual_fcn_prop_gen}. All rate bounds in this comparison that directly analyze ADMM also depend on the dual problem strong convexity and smoothness parameters $\hat{\sigma}$ and $\hat{\beta}$. Therefore, the comparison can be carried out in dual problem parameters $\hat{\sigma}$ and $\hat{\beta}$. Obviously, since our bounds are tight, none of the other bounds can be better. This comparison is merely to show that we improve and/or generalize previous results.

% Due to the duality correspondance between strong convexity and smoothness, we can apply the Douglas-Rachford rates to the dual problem to get corresponding rates for ADMM. 

% Since Douglas-Rachford splitting and ADMM are equivalent in the case
% where $\mathcal{A}=\id$ (that is, Douglas-Rachford is self-dual in the
% sense that
% it gives equivalent algorithms if applied to the primal and the dual
% when $\mathcal{A}=\id$) we can compare Douglas-Rachford convergence rate results
% with ADMM convergence rate results by
% letting $\mathcal{A}=\id$. 
% The convergence rate result in
% \cite{LionsMercier1979} is provided in the
% case of solving monotone inclusion problems using Douglas-Rachford
% splitting, while the other \cite{Davis_Yin_2014,linConvADMM,GhadimiADMM,Panos_acc_DR_2014,laurent_ADMM} consider the convex optimization case
% treated here (which is a subclass of the monotone inclusion problem
% class). 
% We will compare to the results in
% \cite{Davis_Yin_2014,linConvADMM,GhadimiADMM,Panos_acc_DR_2014,laurent_ADMM}
% that hold under Assumption~\ref{ass:prob_ass}. In the monotone
% inclusion problem case considered in \cite{LionsMercier1979},
% Assumption~\ref{ass:prob_ass}
% corresponds to that one of the two maximal monotone operators is
% $\sigma$-strongly monotone and 
% $\beta$-Lipschitz continuous.

In \cite[Proposition~4, Remark~10]{LionsMercier1979}, the linear
convergence rate for Douglas-Rachford splitting with $\alpha=\tfrac{1}{2}$ when solving 
monotone inclusion problems with one operator being $\hat{\sigma}$-strongly
monotone and $\hat{\beta}$-Lipschitz continuous is shown to be
$\sqrt{1-\tfrac{\hat{\sigma}}{2\hat{\beta}}}$. Note that the setting considered in \cite{LionsMercier1979} is more general than the setting in this paper. Recently, this bound was improved in \cite{gisSIAM2015} where tight estimates are provided.

In \cite[Theorem 6]{Davis_Yin_2014}, dual Douglas-Rachford splitting with $\alpha=1$ is shown to converge linearly at least with rate
 $\sqrt{1-\tfrac{\hat{\sigma}}{\hat{\beta}}}$. For $\alpha=\tfrac{1}{2}$ they recover the bound in \cite{LionsMercier1979}. In Figure~\ref{fig:rate_fig}, the former (better) of these rate bounds is plotted. We see that it is more conservative than the one in Corollary~\ref{cor:ADMM_lin_conv_subdiff}.

The optimal convergence rate bound in \cite[Corollary 3.6]{linConvADMM} is
$\sqrt{1/(1+1/\sqrt{\hat{\beta}/\hat{\sigma}})}$ and is analyzed directly in the ADMM setting. Figure~\ref{fig:rate_fig} reveals that our rate bound is tighter.

In \cite{Panos_acc_DR_2014}, the authors show that if the $\gamma$ parameter is small enough
and
if $f$ is a quadratic function, then Douglas-Rachford splitting
is equivalent to a gradient method applied to a smooth convex function named the
Douglas-Rachford envelope. Convergence rate estimates then follow from the gradient method rate estimates. Also, accelerated variants of Douglas-Rachford splitting are proposed, based on fast gradient methods. 
In Figure~\ref{fig:rate_fig}, we compare to the rate bound of the fast Douglas-Rachford splitting in \cite[Theorem 6]{Panos_acc_DR_2014} when applied to the dual. This rate bound is better than the rate bound for standard Douglas-Rachford splitting in \cite[Theorem
4]{Panos_acc_DR_2014}. We note that Corollary~\ref{cor:ADMM_lin_conv_subdiff} gives better rate bounds.

The convergence rate bound provided in \cite{GhadimiADMM} coincides
with the bound
provided in Corollary~\ref{cor:ADMM_lin_conv_subdiff}. The rate bound in \cite{GhadimiADMM} holds for ADMM applied to Euclidean quadratic problems with linear inequality constraints. 
We generalize these results, using a different machinery, to arbitrary real Hilbert spaces (also
infinite dimensional),
to both Douglas-Rachford splitting and ADMM, to general smooth and
strongly convex functions $f$, 
and, perhaps most importantly, to any proper, closed, and convex function $g$.

Finally, we compare our rate bound to the rate bound in \cite{laurent_ADMM}.
Figure~\ref{fig:rate_fig} shows that our bound is tighter. As opposed to all the other
rate bounds in this comparison, the rate bound in \cite{laurent_ADMM} is not
explicit. Rather, a sweep over different rate bound factors is needed.
For each guess, a small semi-definite program is solved
to assess whether the algorithm is guaranteed to converge with that
rate. The quantization level of this sweep is the cause
of the steps in the rate curve in Figure~\ref{fig:rate_fig}.

\section{Metric selection}
\label{sec:metric_select_ADMM}

In this section, we consider problems of the form
\eqref{eq:prob_wA} in an Euclidean setting. We assume $f\in\Gamma_0(\mathbb{H}_M)$, $g\in\Gamma_0(\mathbb{H}_{\hat{K}})$, $\mathcal{A}~:~\mathbb{H}_M\to\mathbb{H}_K$, $\mathcal{B}~:~\mathbb{H}_{\hat{K}}\to\mathbb{H}_K$, $c\in\mathbb{H}_K$ and that:
\begin{ass}~
  \begin{enumerate}[(i)]
  \item $f\in\Gamma_0(\mathbb{H}_M)$ is 1-strongly convex if
    defined on $\mathbb{H}_H$ (i.e., $M=H$) and 1-smooth if defined on $\mathbb{H}_L$ (i.e., $M=L$).
  \item The bounded linear operator
    $\mathcal{A}~:~\mathbb{H}_M\to\mathbb{H}_K$ is surjective.
  \end{enumerate}
  \label{ass:probAss_wA_fd}
\end{ass}
We solve \eqref{eq:prob_wA} by applying Douglas-Rachford splitting on
the dual problem \eqref{eq:dualProb} (or equivalently by applying ADMM
directly on the primal
\eqref{eq:prob_wA}). This algorithm behaves differently depending on which space $\mathbb{H}_K$ the problem is defined and the algorithm is run. We will show how to select a space 
$\mathbb{H}_K$ on which the algorithm rate bound is optimized. To aid in this selection, we show in the
following proposition how the strong 
convexity modulus and smoothness constant of
$d_1\in\Gamma_0(\mathbb{H}_K)$ depend on the space on which it is defined.

\begin{prp}
  Suppose that Assumption~\ref{ass:probAss_wA_fd} holds, that
  $A\in\reals^{m\times n}$ satisfies $Ax=\mathcal{A}x$ for all $x$, and
  that $K=E^TE$. Then
  $d_1\in\mathbb{H}_K$ in \eqref{eq:dFcn} is $\|EAH^{-1}A^TE^T\|_2$-smooth and
  $\lambda_{\min}(EAL^{-1}A^TE^T)$-strongly convex, where
  $\lambda_{\min}(EAL^{-1}A^TE^T)>0$. 
\end{prp}
\begin{pf}
First, we note that $d_1(\mu)=d(\mu)+\langle c,\mu\rangle$ where $d(\mu)=f^*(-\mathcal{A}^*\mu)$ and that a linear function does not change the strong convexity of smoothness parameter of a problem. Therefore, we show the result for $d$.

  First, we relate $\mathcal{A}^*~:~\mathbb{H}_K\to\mathbb{H}_M$ to $A$,
  $M$, and $K$. We have
  \begin{align*}
    \langle\mathcal{A}x,\mu\rangle_K &= \langle Ax,K\mu\rangle_2=\langle
    x,A^TK\mu\rangle_2=\langle x,M^{-1}A^TK\mu\rangle_M\\
    &=\langle M^{-1}A^TK\mu,x\rangle_M = \langle\mathcal{A}^*\mu,x\rangle_M.
  \end{align*}
  Thus, $\mathcal{A}^*\mu = M^{-1}A^TK\mu$ for all $\mu\in\mathbb{H}_K$. 

  Next, we show that the space
  $\mathbb{H}_M$ on which $f$ and $f^*$ are defined does not influence
  the shape of $d$. We denote by $f_H$, $f_L$, and $f_e$ the function
  $f$ defined on $\mathbb{H}_H$, $\mathbb{H}_L$ and $\reals^n$
  respectively and by $\mathcal{A}_H^*~:~\mathbb{H}_K\to\mathbb{H}_H$, $\mathcal{A}_L^*~:~\mathbb{H}_K\to\mathbb{H}_L$, and $A^T~:~\reals^m\to\reals^n$ the
  operator $\mathcal{A}^*$ defined on the different spaces. Further,
  let $d_H:=f_H^*\circ -\mathcal{A}_H^*$,
  $d_L:=f_L^*\circ -\mathcal{A}_L^*$, and $d_e :=
  f_e^*\circ -A^T$. With these definitions both 
  $d_L$ and $d_H$ are defined on $\mathbb{H}_K$, while $d_e$ is defined
  on $\reals^m$. Next we show that $d_L$ and $d_H$ are identical for any
  $\mu$:
  \begin{align*}
    d_H(\mu) &= f_H^*(-\mathcal{A}_H^*\mu) \\&=
    \sup_{x}\left\{\langle-\mathcal{A}_H^*\mu,x\rangle_H-f_H(x)\right\}\\ 
&=\sup_x\left\{\langle-HH^{-1}A^TK\mu,x\rangle_2-f_e(x)\right\}\\
&=\sup_x\left\{\langle-LL^{-1}A^TK\mu,x\rangle_2-f_e(x)\right\}\\
&=\sup_x\left\{\langle-\mathcal{A}_L^*\mu,x\rangle_L-f_L(x)\right\}=d_L(\mu)
  \end{align*}
  where $\mathcal{A}_M^*\mu = M^{-1}A^TK\mu$ is used. This implies that
  we can show properties of $d\in\mathbb{H}_K$ by defining $f$ on any space
  $\mathbb{H}_M$. 
  Thus, Proposition~\ref{prp:dual_fcn_prop_gen} gives that 1-strong
  convexity of $f$ when defined on 
  $\mathbb{H}_H$ implies $\|\mathcal{A}^*\|^2$-smoothness of $d$, where 
  \begin{align*}
    \|\mathcal{A}^*\| &= \sup_{\mu}\left\{\|\mathcal{A}^*\mu\|_{H}~|~\|\mu\|_{K}\leq
      1\right\}\\
    &=\sup_{\mu}\left\{\|H^{-1}A^TK\mu\|_H~|~\|\mu\|_K\leq 1\right\}\\
    &=\sup_{\mu}\left\{\|H^{-1/2}A^TE^TE\mu\|_2~|~\|E\mu\|_2\leq 1\right\}\\
    &=\sup_{\nu}\left\{\|H^{-1/2}A^TE^T\nu\|_2~|~\|\nu\|_2\leq 1\right\}\\
    &=\|H^{-1/2}A^TE^T\|_2.
  \end{align*}
Squaring this gives the smoothness claim.
  To show the strong-convexity claim, we use that 1-smoothness of $f$
  when defined on $\mathbb{H}_L$ implies $\theta^2$-strong convexity of
  $d$ where $\theta>0$ satisfies $\|\mathcal{A}^*\mu\|_L\geq
  \theta\|\mu\|_K$ for all $\mu\in\mathbb{H}_K$, see
  Proposition~\ref{prp:dual_fcn_prop_gen}. We have
  \begin{align*}
    \|\mathcal{A}^*\mu\|_L^2&=\|L^{-1}A^TK\mu\|_L^2\\&=\|L^{-1/2}A^TE^T(E\mu)\|_2^2\\
    &=\|E\mu\|_{EAL^{-1}A^TE^T}^2\\
    &\geq \lambda_{\min}(EAL^{-1}A^TE^T)\|E\mu\|_{2}^2\\
    &= \lambda_{\min}(EAL^{-1}A^TE^T)\|\mu\|_{K}^2,
  \end{align*}
i.e, $\theta^2=\lambda_{\min}(EAL^{-1}A^TE^T)$.
  The smallest eigenvalue $\lambda_{\min}(EAL^{-1}A^TE^T)>0$ since $A$
  is surjective, i.e. has full row rank, and $E$ and $L$ are positive
  definite matrices. This concludes the proof.
\end{pf}

This result shows how the smoothness constant and strong convexity
modulus of $d_1\in\Gamma_0(\mathbb{H}_K)$ change with the space
$\mathbb{H}_K$ on which it is defined. Combining this with
Proposition~\ref{prp:opt_param_subdiff},
we get that the bound on the convergence rate for Douglas-Rachford
splitting applied to the dual problem \eqref{eq:dualProb} (or equivalently ADMM
applied to the primal \eqref{eq:prob_wA}) is optimized by choosing $K=E^TE$ where $E$ is chosen to 
\begin{align}
  \begin{tabular}{ll}
    minimize & $\displaystyle
    \frac{\lambda_{\max}(EAH^{-1}A^TE^T)}{\lambda_{\min}(EAL^{-1}A^TE^T)}$
  \end{tabular}
\label{eq:min_cond_nbr_dual}
\end{align}
and by choosing $\gamma =
\tfrac{1}{\sqrt{\lambda_{\max}(EAH^{-1}A^TE^T)\lambda_{\min}(EAL^{-1}A^TE^T)}}$.
Using a
non-diagonal $K$ usually gives prohibitively expensive
prox-evaluations. Therefore, we propose to select a diagonal $K=E^TE$
that minimizes \eqref{eq:min_cond_nbr_dual}. The reader is referred to
\cite[Section 6]{gisBoydAut2014metric_select} for different 
methods to achieve this exactly and approximately.

% The procedure of selecting the metric for the dual Douglas-Rachford
% splitting is slightly different compared to the primal version. In the
% primal, we assume that we are given matrices $H$ and $L$ that
% define spaces $\mathbb{H}_H$ and $\mathbb{H}_L$ on which the function $f$
% is 1-strongly convex and 1-smooth respectively. To get optimal
% convergence behaviour, we suppose that these matrices are as tight as
% possible. Then a diagonal metric $M$ (that defines $\mathbb{H}_M$ on
% which the problem is defined) is chosen that optimizes the bound on the
% convergence rate. This metric affects only the algorithm, while the
% problem is the same for any choice of metric matrix $M$.

% In the dual formulation, we are not given matrices
% on which the dual is $\beta$-strongly convex and $\sigma$-smooth respectively.
% Tight estimates of these parameters are instead computed and we show
% how these depend on the space $\mathbb{H}_K$ on which the dual
% problem is defined. Since the dual problem actually changes with the
% space on which it is defined, the choice of $\mathbb{H}_K$ affects
% both the shape of the dual problem and the metric used in the
% algorithm. Therefore, the interpretation made for the primal that
% the algorithm metric is chosen to well estimate the level-sets of
% problem is not exactly true in this case. When we select a space
% $\mathbb{H}_K$ for the dual, instead we simultaneously 
% manipulate the level-sets of the dual problem and the metric in the
% algorithm such that the metric well approximates the manipulated level
% sets of the function $d$.

\begin{rem}
It can
be shown (details are omitted for space considerations)  that solving the dual problem on space $\mathbb{H}_K$
using Douglas-Rachford splitting is
equivalent to solving the preconditioned 
problem
\begin{align}
  \begin{tabular}{ll}
    minimize & $f(x)+g(y)$\\
    subject to & $E(Ax+By)=Ec$.
  \end{tabular}
  \label{eq:primPrecond}
\end{align}
using ADMM. Therefore, the metric selection presented here covers the preconditioning suggestion in \cite{GhadimiADMM} as a special case. 
\end{rem}
\begin{rem}
Metric selection and preconditioning is not new for iterative methods. It can be used with the aim to reduce the computational burden within each iteration \cite{ChambollePock,PockChambolle_Diag,Bredies_Sun_2015,Bramble_Uzawa,Hu_nonlin_Uzawa}, and/or with the aim to reduce the total number of iterations like in our analysis and in \cite{GhadimiADMM,gisBoydAut2014metric_select,Benzi_precond,Bramble_Uzawa,Hu_nonlin_Uzawa}. It is interesting to note that in our paper (if we restrict ourselves to quadratic $f$ with Hessian $H$) and in \cite{GhadimiADMM,gisBoydAut2014metric_select,Bramble_Uzawa,Hu_nonlin_Uzawa}, different algorithms are analyzed with different methods, but all analyses suggest to make $AH^{-1}A^T$ well conditioned (by choosing a metric or doing preconditioning). The analyzed algorithms are ADMM here and in \cite{GhadimiADMM}, (fast) dual forward-backward splitting in \cite{gisBoydAut2014metric_select}, and Uzawa's method to solve linear systems of the form $\left[\begin{smallmatrix} H & A^T\\A & 0\end{smallmatrix}\right](x,\mu)=(-q,b)$ in \cite{Bramble_Uzawa,Hu_nonlin_Uzawa}. This linear system is the KKT-condition for the problem $\min_x\{f(x)+g(y)~|~Ax=y\}$ where $f(x)=\tfrac{1}{2}x^THx+q^Tx$ and $g(y)=\iota_{y=b}(y)$. The dual functions are $d_1(\mu) = (q+A^T\mu)^TH^{-1}(A^T\mu+q)$ with Hessian $AH^{-1}A^T$ and $d_2(\mu)=g^*(\mu) = b^T\mu$. The functions $d_1$ in the dual problems in all these analyses are the same (if we restrict ourselves to quadratic $f$ with Hessian $H$ in our setting), but the functions $d_2$ are different. So all the different metric selections (preconditionings) try to make $d_1$, and its Hessian $AH^{-1}A^T$, well conditioned.
\end{rem}

\subsection{Heuristic metric selection}
\label{sec:extensions}

Many problems do not satisfy all assumptions in Assumption~\ref{ass:probAss_wA}, but for many interesting problems a couple of them are typically satisfied. Therefore, we will here discuss metric and parameter selection
heuristics for ADMM (Douglas-Rachford applied to the dual \eqref{eq:dualProb}) when some of the assumptions are not met. We focus on quadratic problems
of the form
\begin{align}
\begin{tabular}{ll}
  minimize&$\displaystyle\underbrace{\tfrac{1}{2}x^TQx+q^Tx+\hat{f}(x)}_{f(x)}+g(y)$\\
subject to & $Ax+By=c$
\end{tabular}
\label{eq:prob_extensions}
\end{align}
where $Q\in\reals^{n\times n}$ is positive semi-definite, $q\in\reals^n$,
$\hat{f}\in\Gamma_0(\reals^n)$, $g\in\Gamma_0(\reals^m)$,
$A\in\reals^{p\times n}$, $B\in\reals^{p\times m}$, and $c\in\reals^p$.

One set of assumptions that guarantee linear convergence for
dual Douglas-Rachford splitting is
that $Q$ is positive definite, that $\hat{f}$ is nonsmooth, and that $A$ has
full row rank. By inverting these, we get a set of assumptions that if anyone of these are satisfied, we cannot guarantee linear convergence:
\begin{enumerate}[(i)]
\item $Q$ is not positive definite, but positive semi-definite.
\item $\hat{f}$ is nonsmooth, e.g., the indicator function of a convex
  constraint set or a piece-wise affine function
\item $A$ does not have full row rank.
\end{enumerate}
In the first case, we loose strong convexity in $f$ and smoothness in the dual $d_1$. In the second case, we loose
smoothness in $f$ and strong convexity $d_1$. In the third case, we loose 
strong convexity in the dual $d_1$.

We directly tackle the case where the assumptions needed to get linear convergence are violated by all points (i), (ii), and (iii). For the dual case we consider, i.e., the ADMM case, we propose to select the
metric as if $\hat{f}\equiv 0$ in \eqref{eq:prob_extensions}. To do this,
we define
the quadratic part of $f$ in \eqref{eq:prob_extensions}
to be $f_{\rm{pc}}(x):=\tfrac{1}{2}x^TQx+q^Tx$ and introduce the function
$d_{\rm{pc}}:=f_{\rm{pc}}^*\circ -A^T$. The conjugate function
of $f_{\rm{pc}}$ is given by
\begin{align*}
  f_{\rm{pc}}^*(y) &= \sup_x\left\{\langle y,x\rangle
    -f_{\rm{pc}}(x)\right\}
  =\begin{cases}
    \tfrac{1}{2}(y-q)^TQ^{\dagger}(y-q) & {\hbox{if
      }}(y-q)\in\mathcal{R}(Q)\\
    \infty & {\hbox{else}}
  \end{cases}
\end{align*}
where $Q^\dagger$ is the pseudo-inverse of $Q$ and $\mathcal{R}$
denotes the range space. This gives
\begin{align*}
  d_{\rm{pc}}(\mu) = \begin{cases}
    \tfrac{1}{2}(A^T\mu+q)^TQ^{\dagger}(A^T\mu-q) & {\hbox{if
      }}(A^T\mu+q)\in\mathcal{R}(Q)\\
    \infty & {\hbox{else}}
  \end{cases}
\end{align*}
with Hessian $AQ^{\dagger}A^T$ on its domain.
We approximate the dual function $d_1(\mu)=f^*(-A^T\mu)+c^T\mu$ with $d_{\rm{pc}}(\mu)+c^T\mu$. Then we propose to select a
diagonal metric
$K=E^TE$ such that this approximation is well conditioned in directions where it has curvature. That is, we propose to
select a metric $K=E^TE$ such that the pseudo condition number of $AQ^\dagger A^T$ is
minimized. This is achieved by finding an $E$ to
\begin{align*}
  \begin{tabular}{ll}
    minimize & $\displaystyle
    \frac{\lambda_{\max}(EAQ^\dagger A^TE^T)}{\lambda_{\min>0}(EAQ^\dagger A^TE^T)}$
  \end{tabular}
\end{align*}
where $\lambda_{\min>0}$ denotes the smallest non-zero eigenvalue. (See \cite[Section
6]{gisBoydAut2014metric_select} for methods to achieve this exactly and approximately.)
This heuristic reduces to the optimal metric choice in the case where linear
convergence is achieved. The $\gamma$-parameter is
also chosen in accordance with the above reasoning and
Corollary~\ref{cor:ADMM_lin_conv_subdiff} as
$\gamma =
\tfrac{1}{\sqrt{\lambda_{\max}(EAQ^\dagger A^TE^T)\lambda_{\min>0}(EAQ^\dagger
  A^TE^T)}}$.

If instead $\hat{f}$ in \eqref{eq:prob_extensions}
is the indicator function of an affine subspace, i.e.,
$\hat{f}=\iota_{Lx=b}$ and $Q$ is strongly convex on the null-space of $L$. Then %$d_1(\mu)=d(\mu)+c^T\mu$ with $d(\mu)$ satisfying
%\begin{align*}
$d_1(\mu) = \tfrac{1}{2}\mu^TAP_{11}A^T\mu+\xi^T\mu+\chi+c^T\mu$
%\end{align*}
where $\xi\in\reals^n$, $\chi\in\reals$, and
%\begin{align}
$\left[\begin{smallmatrix}
Q & L^T\\
L & 0
\end{smallmatrix}\right]^{-1} = \left[\begin{smallmatrix}
P_{11} & P_{12}\\
P_{21} & P_{22}
\end{smallmatrix}\right]$.
%\label{eq:P11_def}
%\end{align}
Then we can choose metric by minimizing the pseudo condition
number of $AP_{11}A^T$ (which is the Hessian of $d_1$) and select
$\gamma$ as \\$\gamma =
\tfrac{1}{\sqrt{\lambda_{\max}(EAP_{11} A^TE^T)\lambda_{\min>0}(EAP_{11}
  A^TE^T)}}$.

\section{Numerical example}

% In this section, we evaluate the proposed metric and parameter selections as well as the heuristic variants. First, we evaluate the theoretically justified method by solving a problem that satisfies the assumptions needed to get linear convergence. Then we evaluate the heuristic method by solving problems arising in a model predictive control setting that do not satisfy the assumptions to get linear convergence.

\subsection{A problem with linear convergence}

\begin{figure}
\centering
%  \psfrag{noprecondboundbound}{\scriptsize{Bound, no precond}}
%  \psfrag{noprecondactualactual}{\scriptsize{Actual, no precond}}
%  \psfrag{precondboundbound}{\scriptsize{Bound, precond}}
%  \psfrag{precondactualactual}{\scriptsize{Actual, precond}}
%  \psfrag{gammas}{$\gamma^\star$}
%  \psfrag{gamma}{$\gamma$}
%  \psfrag{nbritersiters}{$\#$ iters}
\includegraphics{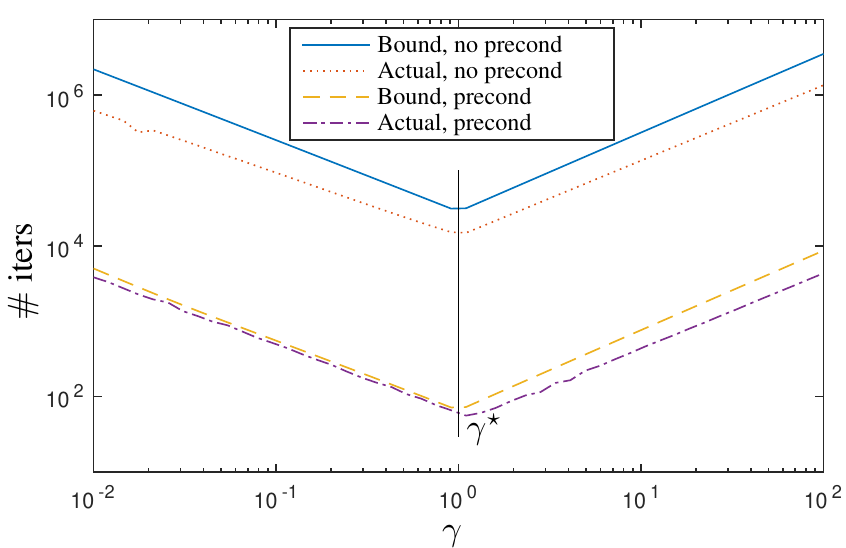}
  \caption{Theoretical bounds on and actual number of iterations to achieve a prespecified relative accuracy of $10^{-5}$ when solving \eqref{eq:weightedL1} using ADMM with $\alpha=1$ for different $\gamma$. We present results with and without preconditioning (metric selection).}
  \label{fig:lin_rate_comparison}
\end{figure}

We consider a weighted Lasso minimization problem formulation with inspiration from \cite{CandesBoydReweightedL1} of the form
\begin{align}
\begin{tabular}{ll}
minimize & $\tfrac{1}{2}\|Ax-b\|_2^2+\|Wx\|_{1}$
\end{tabular}
\label{eq:weightedL1}
\end{align}
with $x\in\reals^{200}$, the data matrix $A\in\reals^{300\times 200}$ is a sparse with on average 10 nonzero entries per row, and $b\in\reals^{300}$. Each nonzero element in $A$ and $b$ is drawn from a Gaussian distribution with zero mean and unit variance. Further, $W\in\reals^{200\times 200}$ is a diagonal matrix with positive diagonal elements drawn from a uniform distribution on the interval $[0,1]$. 

We solve the problem using ADMM in the standard Euclidean setting and in the optimal metric setting according to Section~\ref{sec:metric_select_ADMM}. We use the optimal $\alpha=1$ and different parameters $\gamma$. The $\gamma$ parameter is varied in the range $[10^{-2}\gamma^\star,10^2\gamma^\star]$ where $\gamma^\star$ is the theoretically optimal $\gamma$ in the respective settings (Euclidean and optimal metric). In Figure~\ref{fig:lin_rate_comparison}, the actual number of iterations needed to obtain a specific accuracy (a relative tolerance of $10^{-5}$) as well as the theoretical worst case number of iterations are plotted, both for the Euclidean and optimal metric setting. 

Figure~\ref{fig:lin_rate_comparison} reveals that, for this particular example, the actual numbers of iterations are fairly consistent with the iteration bounds. We also see that there is a lot to gain by running the algorithm with an appropriately chosen metric. Also the optimal parameters $\gamma^\star$ give close to optimal performance in both settings.
% \begin{rem}
% We know from Section~\ref{sec:tigthness_ADMM} that there exist problems for which the theoretical upper bound on the convergence rate coincides with the actual convergence rate. On the other hand, we can also construct problems that always converge in one iteration, provided that $\alpha=1$ and $\gamma$ is chosen optimally. This holds, e.g, for problems where $f=\tfrac{\sigma}{2}\|\cdot\|^2$ with $\sigma>0$, $\mathcal{A}=\id$, $\mathcal{B}=-\id$, $c=0$, and $g$ is an arbitrary proper closed and convex function. This choice always satisfies Assumption~\ref{ass:probAss_wA}. Therefore, the actual convergence for problems within the considered class can vary significantly. 
% \end{rem}

\subsection{A problem without linear convergence}

\begin{figure}
\centering
  \includegraphics[width=0.8\columnwidth]{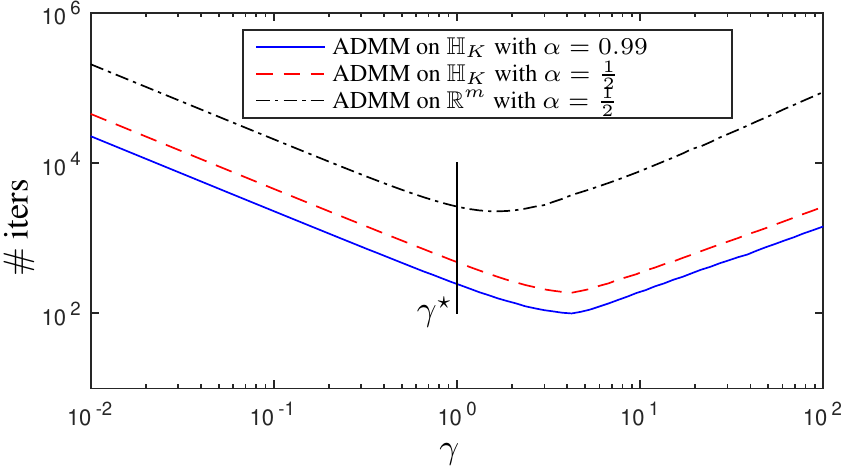}
  \caption{Average number of iterations for different $\gamma$-values,
    different metrics, and different relaxations $\alpha$.}
  \label{fig:gamma_vs_iter}
\end{figure}

Here, we evaluate the heuristic metric and parameter selection method
by applying ADMM to the (small-scale) aircraft control
problem from \cite{Kapasouris,Bemporad97}. As in
\cite{Bemporad97}, the continuous time model from \cite{Kapasouris} 
is sampled using zero-order hold every 0.05~s. The system has four
states $x=(x_1,x_2,x_3,x_4)$, two outputs $y=(y_1,y_2)$, two
inputs $u=(u_1,u_2)$,
and obeys the following  
dynamics
\begin{align*}
  x(t+1) &= \begin{bmatrix}
      0.999  & -3.008 &  -0.113 &  -1.608\\
      -0.000  &  0.986 &   0.048 &   0.000\\
      0.000  &  2.083 &   1.009 &  -0.000\\
      0.000  &  0.053 &   0.050 &   1.000
    \end{bmatrix}x(t)+\begin{bmatrix}
      -0.080  & -0.635\\
      -0.029  & -0.014\\
      -0.868  & -0.092\\
      -0.022  & -0.002
    \end{bmatrix}u(t),\\
  y(t) &= \begin{bmatrix}
      0 & 1 & 0 & 0\\
      0 & 0 & 0 & 1
    \end{bmatrix}x(t).
\end{align*}
The system is unstable, the magnitude of the largest eigenvalue 
of the dynamics matrix is 1.313. The outputs are the
attack and pitch angles, while the inputs are the elevator and
flaperon angles. The inputs are physically constrained to satisfy
$|u_i|\leq 25^\circ$, $i=1,2$. The outputs are soft constrained and
modeled using the piece-wise linear cost function
\begin{align*}
h(y,l,u,s) = \begin{cases} 0 & {\hbox{if }} l\leq y \leq u\\
s(y-u) & {\hbox{if }} y\geq u\\
s(l-y) & {\hbox{if }} y\leq l
\end{cases}
\end{align*}
Especially, the first output is penalized using
$h(y_1,-0.5,0.5,10^6)$ and the second is penalized using
$h(y_2,-100,100,10^6)$. The states and inputs are penalized using
\begin{equation*}
  \ell(x,u,s)=\tfrac{1}{2}\big((x-x_r)^TQ(x-x_r)+u^TRu\big)
\end{equation*}
where $x_r$ is a reference, $Q =
{\rm{diag}}(0,10^2,0,10^2)$, and $R =
10^{-2}I$. Further, the terminal cost is $Q$, and the control and prediction
horizons are $N=10$.
The numerical data in
Figure~\ref{fig:gamma_vs_iter} is obtained by following
a reference trajectory on
the output. The objective is to change the pitch angle from $0^\circ$
to $10^\circ$ and then back to $0^\circ$ while the angle of
attack satisfies the (soft) output constraints $-0.5^\circ \leq y_1 \leq
0.5^\circ$. The constraints on the angle of attack limit the rate on
how fast the pitch angle can be changed. By stacking vectors and
forming appropriate matrices, the full optimization problem
can be written on the form
\begin{align*}
  \begin{tabular}{ll}
    min & $\underbrace{\tfrac{1}{2}z^TQz+
      r_t^Tz+I_{Lz=bx_t}(z)}_{f(z)}+\underbrace{\sum_{i=1}^m
      h(z^{\prime}_i,\underline{d}_i,\bar{d}_i,10^6)}_{g(z^{\prime})}$\\
    s. t. & $Cz=z^\prime$
  \end{tabular}
\end{align*}
where $x_t$ and $r_t$ may change from one sampling instant to the next.
This fits the optimization problem formulation discussed in 
Section~\ref{sec:extensions}. For this problem all items (i)--(iii) violate the assumptions
that guarantee linear convergence.

Since the numerical example treated here
is a model predictive control application, we can spend much
computational effort offline to compute a metric that will be
used in all samples in the online controller. We compute a
diagonal metric $K=E^TE$, where $E$ minimizes the pseudo condition number
of $ECP_{11}C^TE^T$ and $P_{11}$ is implicitly defined by
$\left[\begin{smallmatrix}
Q & L^T\\
L & 0
\end{smallmatrix}\right]^{-1} = \left[\begin{smallmatrix}
P_{11} & P_{12}\\
P_{21} & P_{22}
\end{smallmatrix}\right]$ (as suggested in Section~\ref{sec:extensions}). This matrix $K=E^TE$ defines the space $\mathbb{H}_K$
on which the algorithm is applied. 

In
Figure~\ref{fig:gamma_vs_iter}, the performance
of ADMM when applied on $\mathbb{H}_K$ with relaxations
$\alpha=\tfrac{1}{2}$ and $\alpha=0.99$, and ADMM applied on $\reals^m$
with $\alpha=\tfrac{1}{2}$ is shown. In this
particular example, improvements of about one order of magnitude
are achieved when applied on $\mathbb{H}_K$ compared to when applied
on $\reals^m$.
Figure~\ref{fig:gamma_vs_iter} also shows that ADMM with
over-relaxation performs better than standard ADMM.
The proposed
$\gamma$-parameter selection is denoted by $\gamma^\star$ in
Figure~\ref{fig:gamma_vs_iter} ($E$ or $C$ is scaled to
get $\gamma^\star=1$ for all examples). Figure~\ref{fig:gamma_vs_iter}
shows that $\gamma^\star$ does not coincide with the
empirically found best $\gamma$, but still gives a
reasonable choice of $\gamma$ in all cases.

\section{Conclusions}

We have shown tight global linear convergence rate bounds for
Douglas-Rachford splitting and ADMM. Based on these results, we
have presented methods to select metric and
algorithm parameters for these methods. We
have also provided numerical examples to evaluate the
proposed metric and parameter selection methods for ADMM.

\section{Acknowledgments}

The first author is financially supported by Swedish Foundation for Strategic Research and is a member of the LCCC Linneaus center at Lund University. The anonymous reviewers are gratefully acknowledged for constructive feedback that has led to significant improvements to the manuscript.

\bibliographystyle{plain}
\bibliography{/local/home/pontusg/Research/MPC/papers/references/references}
% \bibliography{references}

\appendix

\section{A lemma}

\label{app:help_lemma}
We state the following lemma that is used in several of the proofs in this appendix.
\begin{lem}
Let $\gamma,\beta,\sigma\in(0,\infty)$ and $\beta\geq\sigma$. Then
\begin{align}
\label{eq:split_max} \max(\tfrac{1-\gamma\sigma}{1+\gamma\sigma},\tfrac{\gamma\beta-1}{1+\gamma\beta})
&=\begin{cases} \tfrac{1-\gamma\sigma}{1+\gamma\sigma} & {\hbox{if }}
  \gamma\leq \tfrac{1}{\sqrt{\beta\sigma}}\\
\tfrac{\gamma\beta-1}{1+\gamma\beta} & {\hbox{if }}
  \gamma\geq \tfrac{1}{\sqrt{\beta\sigma}}
\end{cases}
\end{align}
Further $\max(\tfrac{1-\gamma\sigma}{1+\gamma\sigma},\tfrac{\gamma\beta-1}{1+\gamma\beta})\in[0,1)$
\label{lem:psi_reciprocal}
\end{lem}
\begin{pf}
Let $\psi(x)=\tfrac{1-x}{1+x}$ for $x\geq 0$. This implies $\psi^{\prime}(x)=\tfrac{-2}{(x+1)^2}$, i.e., $\psi$ is (strictly) monotonically deceasing for $x\geq 0$. So for $x\leq 1/y$ and $x\geq 0$, we have
\begin{align*}
\psi(x)=\tfrac{1-x}{1+x}\geq \tfrac{1-1/y}{1+1/y}=-\tfrac{1-y}{y+1}=-\psi(y).
\end{align*}
Similarly if $x\geq 1/y$ and $x\geq 0$, we have
\begin{align*}
\psi(x)=\tfrac{1-x}{1+x}\leq \tfrac{1-1/y}{1+1/y}=-\tfrac{1-y}{y+1}=-\psi(y).
\end{align*}
So for $x\geq 0$, $\psi(x)\leq -\psi(y)$ if and only if $x\geq 1/y$. Therefore
\begin{align}
\nonumber \max(\tfrac{1-\gamma\sigma}{1+\gamma\sigma},\tfrac{\gamma\beta-1}{1+\gamma\beta})
&= \max(\psi(\gamma\sigma)),-\psi(\gamma\beta))\\
\nonumber &=\begin{cases} \psi(\gamma\sigma) & {\hbox{if }}
  0< \gamma\leq \tfrac{1}{\sqrt{\beta\sigma}}\\
-\psi(\gamma\beta) & {\hbox{if }}
  \gamma\geq \tfrac{1}{\sqrt{\beta\sigma}}
\end{cases}\\
\label{eq:split_max}&=\begin{cases} \tfrac{1-\gamma\sigma}{1+\gamma\sigma} & {\hbox{if }}
  0< \gamma\leq \tfrac{1}{\sqrt{\beta\sigma}}\\
\tfrac{\gamma\beta-1}{1+\gamma\beta} & {\hbox{if }}
  \gamma\geq \tfrac{1}{\sqrt{\beta\sigma}}
\end{cases}
\end{align}
Finally, since $\psi(0)=1$ and $\psi(1)=0$, and $\psi$ is strictly monotonically decreasing, $\psi(x)\in[0,1)$ for $x\in(0,1]$. Also $\psi(x)\to -1$ as $x\to\infty$, so $\psi(x)\in(-1,0]$ for $x\geq 1$.  Therefore, since $\beta\geq\sigma$, $\psi(\gamma\sigma)\in[0,1)$ if $0\leq\gamma\leq\tfrac{1}{\sigma}$ and $-\psi(\gamma\beta)\in[0,1)$ if $\gamma\geq\tfrac{1}{\beta}$. Since $\tfrac{1}{\beta}\leq\tfrac{1}{\sigma}$, \eqref{eq:split_max} implies that $\max(\tfrac{1-\gamma\sigma}{1+\gamma\sigma},\tfrac{\gamma\beta-1}{1+\gamma\beta})=\max(\psi(\gamma\sigma),-\psi(\gamma\beta))\in[0,1)$.
\end{pf}

\section{Proof to Proposition~\ref{prp:prox_coco}}

\label{app:prox_coco_pf}
\begin{pf}
Since $f$ is $\sigma$-strongly convex and $\beta$-smooth and since $\gamma\in(0,\infty)$, $f_{\gamma}$
is $(1+\gamma\sigma)$-strongly convex and $(1+\gamma\beta)$-smooth.
Therefore \cite[Theorem~18.15]{bauschkeCVXanal} and
\cite[Theorem~13.32]{bauschkeCVXanal} directly imply that
$f_{\gamma}^*$ is $\tfrac{1}{1+\gamma\sigma}$-smooth and
$\tfrac{1}{1+\gamma\beta}$-strongly convex. From the 
smoothness definition in Definition~\ref{def:smoothness}, we get that
\begin{align}
\label{eq:smooth_exp}\tfrac{1}{2(1+\gamma\sigma)}&\|\cdot\|^2-f_{\gamma}^*
=\left(\tfrac{1}{2(1+\gamma\sigma)}-\tfrac{1}{2(1+\gamma\beta)}\right)\|\cdot\|^2-(f_{\gamma}^*-\tfrac{1}{2(1+\gamma\beta)}\|\cdot\|^2)
\end{align}
is convex. Further, Definition~\ref{def:strConv} implies that
$f_{\gamma}^*-\tfrac{1}{2(1+\gamma\beta)}\|\cdot\|^2$ is convex (since $f_{\gamma}^*$ is $\tfrac{1}{1+\gamma\beta}$-strongly convex), and therefore
  \eqref{eq:smooth_exp} is the definition of
  $\tfrac{1}{2(1+\gamma\sigma)}-\tfrac{1}{2(1+\gamma\beta)}$-smoothness of $f_{\gamma}^*-\tfrac{1}{2(1+\gamma\beta)}\|\cdot\|^2$. 

Now, let $\beta=\sigma$, then $f_{\gamma}^*-\tfrac{1}{2(1+\gamma\beta)}\|\cdot\|^2$ is 0-smooth, or equivalently by applying \cite[Theorem~18.15]{bauschkeCVXanal}, $\nabla f_{\gamma}^*-\tfrac{1}{1+\gamma\beta}\id={\rm{prox}}_{\gamma f}-\tfrac{1}{1+\gamma\beta}\id$ is 0-Lipschitz. Let $\beta>\sigma$, then \cite[Theorem~18.15]{bauschkeCVXanal} implies that $\tfrac{1}{2(1+\gamma\sigma)}-\tfrac{1}{2(1+\gamma\beta)}$-smoothness of $f_{\gamma}^*-\tfrac{1}{2(1+\gamma\beta)}\|\cdot\|^2$ is equivalent to $\tfrac{1}{\tfrac{1}{1+\gamma\sigma}-\tfrac{1}{1+\gamma\beta}}$-cocoercivity of $\nabla f_{\gamma}^*-\tfrac{1}{1+\gamma\beta}\id={\rm{prox}}_{\gamma f}-\tfrac{1}{1+\gamma\beta}\id$. This concludes the proof.
\end{pf}

\section{Proof to Theorem~\ref{thm:refl_res_contr_subdiff}}

\label{app:refl_res_contr_pf}

\begin{pf}
First assume that $\beta>\sigma$. By Proposition~\ref{prp:prox_coco}, ${\rm{prox}}_{\gamma f}-\tfrac{1}{1+\gamma\beta}\id$ is $\tfrac{1}{\tfrac{1}{1+\gamma\sigma}-\tfrac{1}{1+\gamma\beta}}$-cocoercive. Therefore, according to \eqref{eq:coco_expression}:
\begin{align}
{\rm{prox}}_{\gamma f}-\tfrac{1}{1+\gamma\beta}\id=\tfrac{1}{2}(\tfrac{1}{1+\gamma\sigma}-\tfrac{1}{1+\gamma\beta})(\id+C)
\label{eq:prox_characterization}
\end{align}
for some nonexpansive mapping $C$. 

When $\beta=\sigma$, Proposition~\ref{prp:prox_coco} implies that ${\rm{prox}}_{\gamma f}=\tfrac{1}{1+\gamma\beta}\id$. Therefore, also in this case, it can be represented by \eqref{eq:prox_characterization} (since the right hand side is 0).
We get
\begin{align*}
\|R_{\gamma f}x-R_{\gamma f}y\|&=\|(2{\rm{prox}}_{\gamma f}-\id)x-(2{\rm{prox}}_{\gamma g}-\id)y\|\\
&=\|(\tfrac{1}{1+\gamma\sigma}-\tfrac{1}{1+\gamma\beta}+\tfrac{2}{1+\gamma\beta}-1)(x-y)\\
&\quad+(\tfrac{1}{1+\gamma\sigma}-\tfrac{1}{1+\gamma\beta})(Cx-Cy)\|\\
&=\|(\tfrac{1}{1+\gamma\sigma}+\tfrac{1}{1+\gamma\beta}-1)(x-y)\\
&\quad+(\tfrac{1}{1+\gamma\sigma}-\tfrac{1}{1+\gamma\beta})(Cx-Cy)\|\\
&\leq(|\tfrac{1}{1+\gamma\sigma}+\tfrac{1}{1+\gamma\beta}-1|+\tfrac{1}{1+\gamma\sigma}-\tfrac{1}{1+\gamma\beta})\|x-y\|.
\end{align*}
So $R_{\gamma f}$ is Lipschitz continuous with constant $|\tfrac{1}{1+\gamma\sigma}+\tfrac{1}{1+\gamma\beta}-1|+\tfrac{1}{1+\gamma\sigma}-\tfrac{1}{1+\gamma\beta}$.
The kink in the absolute value term is when
\begin{align*}
&&0&=\tfrac{1}{1+\gamma\sigma}+\tfrac{1}{1+\gamma\beta}-1\\
\Leftrightarrow&& 0&=1+\gamma\beta+1+\gamma\sigma-(1+\gamma\sigma)(1+\gamma\beta)\\
\Leftrightarrow&& 0&=1-\gamma^2\sigma\beta,
\end{align*}
i.e, when $\gamma=\tfrac{1}{\sqrt{\beta\sigma}}$. For $\gamma\in(0,\tfrac{1}{\sqrt{\beta\sigma}}]$, the expression in the absolute value is positive, and the Lipschitz constant is
\begin{align*}
\tfrac{1}{1+\gamma\sigma}+\tfrac{1}{1+\gamma\beta}-1+\tfrac{1}{1+\gamma\sigma}-\tfrac{1}{1+\gamma\beta}
=\tfrac{2}{1+\gamma\sigma}-1=\tfrac{1-\gamma\sigma}{1+\gamma\sigma}.
\end{align*}
For $\gamma\in[\tfrac{1}{\sqrt{\beta\sigma}},\infty)$, the expression in the absolute value is negative, and the Lipschitz constant is
\begin{align*}
1-\tfrac{1}{1+\gamma\sigma}-\tfrac{1}{1+\gamma\beta}+\tfrac{1}{1+\gamma\sigma}-\tfrac{1}{1+\gamma\beta}
=1-\tfrac{2}{1+\gamma\beta}=\tfrac{\gamma\beta-1}{1+\gamma\beta}.
\end{align*}
Lemma~\ref{lem:psi_reciprocal} in Appendix~\ref{app:help_lemma} implies that the Lipschitz constant therefore can be written as
\begin{align*}
\max(\tfrac{1-\gamma\sigma}{1+\gamma\sigma},\tfrac{\gamma\beta-1}{\gamma\beta+1})
\end{align*}
and that $\max(\tfrac{1-\gamma\sigma}{1+\gamma\sigma},\tfrac{\gamma\beta-1}{\gamma\beta+1})\in[0,1)$, i.e. that $R_{\gamma f}$ is a contraction.
This concludes the proof.

\end{pf}

  \section{Proof to Theorem~\ref{thm:DR_lin_conv_subdiff}}

  \label{app:DR_lin_conv_subdiff_pf}

  \begin{pf} 
    Since $\gamma\in(0,\infty)$, \cite[Corollary 23.10]{bauschkeCVXanal} implies that $R_{\gamma g}$ is
    nonexpansive and by Theorem~\ref{thm:refl_res_contr_subdiff},
    $R_{\gamma f}$ is 
    $\delta=\max(\tfrac{1-\gamma\sigma}{1+\gamma\sigma},\tfrac{\gamma\beta-1}{\gamma\beta+1})$-contractive. Therefore
    the composition $R_{\gamma g}R_{\gamma f}$ is also
    $\delta$-contractive since
    \begin{align}
      \|R_{\gamma g}R_{\gamma f}z_1-R_{\gamma g}R_{\gamma
        f}z_2\|&\leq\|R_{\gamma f}z_1-R_{\gamma f}z_2\|
      \leq\delta\|z_1-z_2\|
      \label{eq:reflResolvCompContr}
    \end{align}
    for any $z_1,z_2\in\hilbert$.
    Now, let $T = (1-\alpha)I+\alpha R_{\gamma g}R_{\gamma f}$ be the
    Douglas-Rachford
    operator in \eqref{eq:genDougRach}. Since $\bar{z}$ is a fixed-point
    to $R_{\gamma g}R_{\gamma f}$ it is also a fixed-point to $T$, i.e.,
    $\bar{z}=T\bar{z}$. Thus
    \begin{align*}
      \|z^{k+1}-\bar{z}\|&=\|Tz^{k}-T\bar{z}\|^2\\
      &=\|(1-\alpha)(z^k-\bar{z})+\alpha (R_{\gamma g}R_{\gamma
        f}z^k-R_{\gamma g}R_{\gamma f}\bar{z})\|\\
      &\leq|1-\alpha|\|z^k-\bar{z}\|+\alpha
      \|R_{\gamma g}R_{\gamma f}z^k-R_{\gamma g}R_{\gamma f}\bar{z}\|\\
      &\leq\left(|1-\alpha|+\alpha\delta\right)\|z^k-\bar{z}\|\\
      &=\left(|1-\alpha|+\alpha\max(\tfrac{1-\gamma\sigma}{1+\gamma\sigma},\tfrac{\gamma\beta-1}{\gamma\beta+1})\right)\|z^k-\bar{z}\|
    \end{align*}
    where \eqref{eq:reflResolvCompContr} is used in the second inequality. For any $\gamma\in(0,\infty)$, Lemma~\ref{lem:psi_reciprocal} says that $\delta=\max(\tfrac{1-\gamma\sigma}{1+\gamma\sigma},\tfrac{\gamma\beta-1}{\gamma\beta+1})\in[0,1)$ and straightforward manipulations show that the factor 
\begin{align*}
|1-\alpha|+\alpha\delta<1
\end{align*}
if and only if $\alpha\in (0,\tfrac{2}{1+\delta})$. This concludes the proof.
  \end{pf}

\section{Proof of Corollary~\ref{cor:x_lin_conv}}

\label{app:x_lin_conv_pf}

\begin{pf}
Since $f$ is $\sigma$-strongly convex, $f_{\gamma}=\gamma f+\tfrac{1}{2}\|\cdot\|^2$ is $(1+\gamma\sigma)$-strongly convex, and $\nabla f_{\gamma}^*$ is $\tfrac{1}{1+\gamma\sigma}$-Lipschitz continuous, see \cite[Proposition~18.15]{bauschkeCVXanal}. 

Now, recall from Proposition~\ref{prp:prox_to_conj} that $\nabla f_{\gamma}^*(z)={\rm{prox}}_{\gamma f}(z)$. Further, let $\bar{z}$ be a fixed-point to $R_{\gamma g}R_{\gamma f}$, i.e., $\bar{z}$ satisfies $\bar{z}=R_{\gamma g}R_{\gamma f}\bar{z}$ and $x^\star={\rm{prox}}_{\gamma f}(\bar{z})$, see \eqref{eq:DR_optcond}. Therefore
\begin{align*}
\|x^{k+1}-x^\star\|&=\|{\rm{prox}}_{\gamma f}(z^{k+1})-{\rm{prox}}_{\gamma f}(\bar{z})\|\\
&\leq \tfrac{1}{1+\gamma\sigma}\|z^{k+1}-\bar{z}\|\\
&\leq \tfrac{1}{1+\gamma\sigma}(|1-\alpha|+\delta\alpha)\|z^k-\bar{z}\|\\
&\leq \tfrac{1}{1+\gamma\sigma}(|1-\alpha|+\delta\alpha)^{k+1}\|z^{0}-\bar{z}\|.
\end{align*}
where the second and third inequalities come from Theorem~\ref{thm:DR_lin_conv_subdiff}. This concludes the proof.
\end{pf}

\section{Proof of Theorem~\ref{thm:tight_bound}}
\label{app:tight_bound_pf}

Let $z^0=(1,0)$ or $z^0=(0,1)$. Then separability of $R_{\gamma f}$ \eqref{eq:f_refl_prox} and $R_{\gamma g_1}=\id$ implies that the Douglas-Rachford algorithm \eqref{eq:genDougRach} for minimizing $f+g_1$ is
\begin{align}
z^{k+1}=\left((1-\alpha)+\alpha\tfrac{1-\gamma\lambda_i}{1+\gamma\lambda_i}\right)z^0
\label{eq:DR_g1}
\end{align}
with $\lambda_i=\beta$ if $z^0=(1,0)$ and $\lambda_i=\sigma$ if $z^0=(0,1)$. The exact rate is 
\begin{align}
\left|1-\alpha+\alpha\tfrac{1-\gamma\lambda_i}{1+\gamma\lambda_i}\right|.
\label{eq:DR_g1_rate}
\end{align}
Similarly, when minimizing $f+g_2$ we get (since $R_{\gamma g_2}=-\id$ and $R_{\gamma f}$ is linear)
\begin{align}
z^{k+1}=\left((1-\alpha)-\alpha\tfrac{1-\gamma\lambda_i}{1+\gamma\lambda_i}\right)z^0
\label{eq:DR_g2}
\end{align}
with exact rate
\begin{align}
\left|1-\alpha-\alpha\tfrac{1-\gamma\lambda_i}{1+\gamma\lambda_i}\right|.
\label{eq:DR_g2_rate}
\end{align}

Let's consider the following four cases
\begin{enumerate}[(i)]
\item $\alpha\leq 1$, $\gamma\in(0,\tfrac{1}{\sqrt{\beta\sigma}}]$, $g=g_1$, $z_0=(0,1)$ $\Rightarrow$ $\lambda_i=\sigma$
\item $\alpha\geq 1$, $\gamma\geq\tfrac{1}{\sqrt{\beta\sigma}}$, $g=g_1$, $z_0=(1,0)$ $\Rightarrow$ $\lambda_i=\beta$
\item $\alpha\leq 1$, $\gamma\geq\tfrac{1}{\sqrt{\beta\sigma}}$, $g=g_2$, $z_0=(1,0)$ $\Rightarrow$ $\lambda_i=\beta$
\item $\alpha\geq 1$, $\gamma\in(0,\tfrac{1}{\sqrt{\beta\sigma}}]$, $g=g_2$, $z_0=(0,1)$ $\Rightarrow$ $\lambda_i=\sigma$
\end{enumerate}
where $\lambda_i$ refers to the $\lambda_i$ in \eqref{eq:DR_g1}-\eqref{eq:DR_g2_rate}.

Using Lemma~\ref{lem:psi_reciprocal} in Appendix~\ref{app:help_lemma}, it is straightforward to verify that the exact rate factors \eqref{eq:DR_g1_rate} for (i) and (ii) and \eqref{eq:DR_g2_rate} for (iii) and (iv) equals the upper bound on the rate in Theorem~\ref{thm:DR_lin_conv_subdiff} for $\alpha\in(0,\tfrac{2}{1+\delta})$ with $\delta$ in \eqref{eq:delta}. 

It is also straightforward to verify that for $\alpha\not\in(0,\tfrac{2}{1+\delta})$, the exact rate is greater than or equal to one, so it does not converge.

We show this for the first case, and leave the three other very similar cases for the reader to verify. For case (i) with $\alpha\in(0,1]$, Lemma~\ref{lem:psi_reciprocal} implies that the rate \eqref{eq:DR_g1_rate} is
\begin{align*}
\left|1-\alpha+\alpha\tfrac{1-\gamma\sigma}{1+\gamma\sigma}\right|&=1-\alpha+\alpha\max(\tfrac{1-\gamma\sigma}{1+\gamma\sigma},\tfrac{\gamma\beta-1}{1+\gamma\beta})\\
&=|1-\alpha|+\alpha\max(\tfrac{1-\gamma\sigma}{1+\gamma\sigma},\tfrac{\gamma\beta-1}{1+\gamma\beta})
\end{align*}
which is the optimal rate in Theorem~\ref{thm:DR_lin_conv_subdiff}. If instead $\alpha\leq 0$, we get
\begin{align*}
\left|1-\alpha+\alpha\tfrac{1-\gamma\sigma}{1+\gamma\sigma}\right|&=1+|\alpha|(1-\max(\tfrac{1-\gamma\sigma}{1+\gamma\sigma},\tfrac{\gamma\beta-1}{1+\gamma\beta}))\geq 1
\end{align*}
since $\max(\tfrac{1-\gamma\sigma}{1+\gamma\sigma},\tfrac{\gamma\beta-1}{1+\gamma\beta})\in[0,1)$ by Lemma~\ref{lem:psi_reciprocal}. So it does not converge to the solution.

Repeating similar arguments for the three other cases gives the result.

\vfill
\end{document}